\documentclass[preprint,12pt]{elsarticle}

\usepackage{graphicx}														\usepackage{epstopdf}				
\usepackage{amsmath}
\usepackage{amssymb}
\usepackage{array}
\usepackage[caption=false,font=normalsize,labelfont=sf,textfont=sf]{subfig}
\usepackage{stfloats}

\newtheorem{theorem}{Theorem}
\newtheorem{lemma}{Lemma}
\newtheorem{definition}{Definition}
\newtheorem{assumption}{Assumption}
\newtheorem{algorithm}{Algorithm}
\newtheorem{remark}{Remark}

\journal{System $\&$ Control Letters}

\begin{document}

\begin{frontmatter}

\title{An Output-Feedback Control Approach to the $H_{\infty}$ Consensus Integrated with Transient Performance Improvement Problem}

\author[label1]{Jingyao Wang}
\ead{wangjingyao1@xmu.edu.cn}
 \author[label1,label2]{Zhisheng Duan}
 \ead{duanzs@pku.edu.cn}
 \author[label1]{Jianping Zeng}
 \ead{jpzeng@xmu.edu.cn}
 \address[label1]{Department of Automation, Xiamen University, Xiamen 361000, China}
 \address[label2]{Department of Technology and Engineering Science,
Peking University, Beijing 100871, China
}

\begin{abstract}
This paper considers the consensus performance improvement problem of networked general linear agents subject to external disturbances over Markovian randomly switching communication topologies. The consensus control laws can only use its local output information. Firstly, a class of full-order observer-based control protocols is proposed to solve this problem, which depends solely on the relative outputs of neighbours. Then, to eliminate the redundancy involved in the full-order observer, a class of reduced-order observer-based control protocols is designed. Algorithms to construct both protocols are presented, which guarantee that agents can reach consensus in the asymptotic mean square sense when they are not perturbed by disturbances, and that they have decent $H_{\infty}$ performance and transient performance when the disturbances exist.
At the end of this manuscript, numerical simulations which apply both algorithms to four networked Raptor-$90$ helicopters are performed to verify the theoretical results.

\end{abstract}


\end{frontmatter}

\footnotetext{This work was partially supported by the National Nature Science Foundation of China under Grant Nos. 61803319, U1713223 and 61673026.}

\section{Introduction}\label{sec_introduction}
The coordination control of networked agents, for example, unmanned aerial vehicles (UAVs), has gained increasing attention thanks to its paramount importance in both civilian applications like monitoring oil fields and pipelines, and applications to homeland security like border patrol. Due to this, numerous control protocols have been constructed to address various control problems, e.g., \cite{Liu2015finite,Zhou1996robust}. Though these works have provided a theoretical underpinning to the application of coordination control, there is still a critical request for coordination control laws which are uncomplicated to implement in practice and designed for hazardous environments.

In practice, agents, for example, UAVs, when they are required to conduct tasks like monitoring oil fields and pipelines, are usually exposed to complicated environments. Under such circumstance, the wind in the outdoor space and other sources of external disturbances may deteriorate the control performance and even threaten the life of agents. For example, small fixed-wing UAVs are highly susceptible to wind. Moreover, agents are also suffering from disturbances cased by issues like un-modeled high-order dynamics, which arise quite often in practice and are to be found, for instance, in aircraft control systems and large scale flexible structures for space stations (such as antenna, solar arrays, etc.). Thus, guaranteeing good robustness performance of agents is highly valued in the design of control laws. To improve the robustness performance of agents, great efforts have been devoted and hence there are a large amount of references on this topic in the published literature. According to the number of agents, results on this topic can be classified into robust control design for single agent and robust control design for a group of agents. For the single agent case, the novel publication \cite{Zhou1996robust} provides benchmark solutions, like $H_{\infty}$ optimal control and $H_2$ optimal control, to deal with various external disturbance attenuation problems. For the networked agents case, 
there are also many researchers who have endeavored to improve the $H_{\infty}$ and $H_2$ consensus performances of multi-agent systems with linear dynamics. To name a few, \cite{robust_H_infty_containment_Jia_2015,wang2013h,robust_consensus_Zeng_2015,Hinf_Htwo_perfor_region_2010,
Wang2015distributed,Wang2015consensus,lin2017distributedconstrained,zhao2012hinfinity}. 
Note that a common assumption is adopted in the above references, i.e., the communication graph is time-invariant or fixed and even known previously. However, this assumption has limited the application of the proposed control algorithms in \cite{robust_H_infty_containment_Jia_2015,wang2013h,robust_consensus_Zeng_2015,Hinf_Htwo_perfor_region_2010,Wang2015distributed,Wang2015consensus,lin2017distributedconstrained,zhao2012hinfinity}. In fact, the communication topologies are usually randomly switching, because of many practical backgrounds including gossip algorithms and communication patterns (for example, \cite{Boyd2006randomized,Matei2009consensus_problems}). In addition, during the information transmission, packet drop and node failure phenomena can be described as random switching graph process. Furthermore, the importance of the multi-agent consensus with various random graph has attracted a lot of interest recently. Results related to this include but not limited to \cite{Li2015containment,Kan2016leader,Wang2015seeking,you2013consensus,Pan2017distributedcooperative}.

To the best of our knowledge, there are only a few works investigating the robustness performance of multi-agent systems with the randomly switching communication topologies. In \cite{Mu2014L}, the $L_2$-$L_{\infty}$ containment control problem of second-order multi-agent systems with Markovian switching topologies and external disturbances is investigated.
In \cite{Wang2017distributed}, we proposed a class of state-feedback control protocols to improve the $H_{\infty}$ and transient performances of agents with linear dynamics over randomly switching topologies. We mainly focused on improving the transient performance of agents in networks, i.e., reducing the big overshoot and large oscillation, which are caused by the initial-state uncertainty. The reason of improving the transient performance is that, as pointed out in \cite{IJC_2010}, because of the initial-state uncertainty, the transient performance might be unacceptable or might deteriorate the disturbance rejection level of $H_{\infty}$ control. Following this pattern, we also studied the target containment control problem of multi-agent systems under the Markovian randomly switching topologies in \cite{wang2018h}. It is worth noting that the control laws in \cite{Mu2014L, Wang2017distributed,wang2018h} rely on the relative states of neighbouring agents. In fact, the state information might be unavailable in many circumstances. Instead, each agent can access the local output information, i.e., the output of itself and the outputs from its neighbours. Moreover, the relative outputs of neighbours can be measured by sensor networks, i.e., a variety of sonar sensors. Therefore, designing control laws which solely depend on local output information is more practical.

Motivated by the aforementioned works, we study the $H_{\infty}$ consensus control integrated with transient performance problem of general linear multi-agent systems with Markovian randomly switching interconnections. Here, we first propose a class of distributed control protocols which uses an observer to estimate the entire state information.  The main advantage of this class of control protocols is that it
solely requires the relative output information of neighbouring agents. 
However, the main drawback of these protocols is that estimating the entire state information
may involve a certain degree of redundancy, since actually part of the state information is available in the output. To eliminate the redundancy, we design a class of reduced-order observer to estimate the missing state information and then design the consensus control laws based on the relative states of observers and the relative outputs of neighbouring agents. The main difficulty of constructing the above two classes of control laws lies in that each possible communication graph is directed, which implies that the Laplacian matrices of these graphs are generally asymmetric, rendering the constructing of the consensus protocols and the selection of appropriate Lyapunov function far from being easy.

The reminder of this paper is as follows. Section II presents some preliminary results which are key to solve the concerned problems. Sections III and IV are dedicated to address the $H_{\infty}$ stochastic consensus integrated with transient performance problem. Section V shows the effectiveness of the proposed control laws with a simulation example.

\textbf{Notation:}
let $\mathrm{R}^{n\times n}$ be the set of $n\times n$ real matrices. Let $\mathbb{L}^{r}_2[0,\infty)$ denote the space of square integrable vector functions over $[0,\infty)$, which are of $r$ dimension. And let $I_N$ denote the identity matrix of order $N$. Represent the $L_2$ norm of the corresponding function with $\|\cdot \|$.

\section{Preliminaries}\label{sec_notation}

Denote by $\sigma(t),\ t\in\mathrm{R}_{+}$, a right-continuous time-homogeneous Markov process based on the probability space. $\sigma(t)$ takes values in a finite state space $S\! = \! \{ 1,2,\ldots,s\}$. Its infinitesimal generator matrix is denoted by $Q\! = \! [q_{ij}]$, where
$q_{ij}$ denotes the transition rate from state $i$ to state $j$ with $q_{ij}\! \geq \! 0$ if $i\! \neq \! j$; $q_{ii}\! = \! \sum\limits_{j\! \neq \! i}^N q_{ij}$ otherwise. Thus, matrix $Q$ is a transition matrix with row summation being zero and all off-diagonal elements non-negative.

At time $t$, denote the communication graph among all agents (nodes) by $\mathcal{G}_{\sigma(t)}\! \triangleq \! (\mathcal{V},\ \mathcal{E}_{\sigma(t)})$, where $\mathcal{V}\! \triangleq \! \{1,\ 2,\ldots,\ N \}$ and $\mathcal{E}_{\sigma(t)}$ are respectively the node set and the edge set. The associated adjacency matrix is $\mathcal{A}(\sigma(t))\! \triangleq \! [a_{ij}(\sigma(t))]\in\mathrm{R}^{N\! \times \! N}$, where $a_{ij}(\sigma(t))\! > \! 0$ if $(j,\ i)\in\mathcal{E}_{\sigma(t)}$; $a_{ij}(\sigma(t))\! = \! 0$ otherwise. The Laplacian matrix associated with $\mathcal{G}_{\sigma(t)}$ can thus be defined, i.e.,  $\mathcal{L}_{\sigma(t)}\! \triangleq \! [l_{ij}({\sigma(t)})]\! \in \! \mathrm{R}^{N\! \times \! N}$, where $l_{ij}({\sigma(t)})\! = \! - \! a_{ij}({\sigma(t)})$ if $i\! \neq \! j$; $l_{ii}({\sigma(t)})\! = \! \sum_{j\! \neq \! i}^N a_{ij}({\sigma(t)})$ otherwise. 
We say the communication graph at time $t$ is balanced, if each node's in-degree equals to its out-degree.

Represent the union of the $s$ graphs $\mathcal{G}_i\! = \! (\mathcal{V},\mathcal{E}_i),\ i=1,2,\ldots,s, $ by $\mathcal{G}_{\text{un}}\triangleq \cup_{i=1}^s \mathcal{G}_i\! =\! (\mathcal{V},\ \cup_{i=1}^s\mathcal{E}_i)$. The Laplacian matrix associated with the union graph $\mathcal{G}_{\text{un}}$ is denoted by $\mathcal{L}_{\text{un}}$, which is the sum of $\mathcal{L}_i,\ i=1,2,\ldots,s$.

\begin{lemma}\label{lemma_Yang}
For any vectors $p,\ q\in \mathrm{R}^n$ and positive definite matrix $\Phi\in\mathrm{R}^{n \times n}$, the following matrix inequality holds:
\begin{equation}
2p^Tq\leq p^T\Phi p + q^T\Phi^{-1}q.
\end{equation}
\end{lemma}

\section{Problem Formulation and Full-Order Observer-Based Control Protocol}
This section focuses on formulating the $H_{\infty}$ consensus control integrated with transient performance improvement problem and then constructing the distributed control protocol to deal with it.

\subsection{Problem Formulation}\label{sec_main_one_one}
Consider the following $N$ agents:
\begin{eqnarray}\label{dynamics_each_agent}
\dot{x}_i(t)&\! = \! & Ax_i(t) \! + \! B u_i(t)\! + \! D \omega_i(t),\nonumber\\
y_i(t)&\! = \! & C_1 x_i(t), \ i=1,2,\ldots,N,
\end{eqnarray}
where $x_i(t)\in \mathrm{R}^n,\ u_i(t)\in \mathrm{R}^m, \omega_i(t)\in \mathbb{L}_2^\ell[0,\ \infty]$ and $y_i(t)\in \mathrm{R}^{q_1} $ are the state, the control input, the external disturbance and the measured output of the $i$-th agent, respectively;
$A\in \mathrm{R}^{n \times n},\ B \in \mathrm{R}^{n \times m},\ C_1 \in \mathrm{R}^{q_1 \times n}$ and $D \in \mathrm{R}^{n \times \ell}$ are constant matrices.

In this paper, we consider the network of agents described by (\ref{dynamics_each_agent}) with Markovian switching communication graph 
 $\mathcal{G}_{\sigma(t)}$. The switching process is defined in Section \ref{sec_notation}  and satisfies the following assumption:
\begin{assumption}
The continuous-time Markov process with a transition rate matrix $Q$ governing the switching process of the communication topologies is ergoic.
\end{assumption}

As a result of Assumption $1$, for the time-homogeneous Markov process, its invariant distribution is unique $\pi=[\begin{array}{*{20}{c}}\pi_1 & \pi_2 & \ldots \pi_s \end{array}]^T$, where $\pi_i\geq 0$, for all $i\in\mathcal{S}$; thus the distribution of $\sigma(t) $ equals to $\pi$ for all $t\in\mathrm{R}_{+}$.

Herein the union $\mathcal{G}_{\text{un}}$ of the $s$ graphs satisfies the following assumption:
\begin{assumption}
There exists a directed spanning tree in the union graph $\mathcal{G}_{\text{un}}$ of all the possible topologies, and each possible topology is directed and balanced.
\end{assumption}

Now we proceed to present the concerned problem.

Considering the fact that only local outputs are available, we introduce the following observer for each agent to estimate its state information:
\begin{eqnarray}\label{FS_observer_each_agent}
\dot{\hat{x}}_i(t)&\! = \! & A \hat{x}_i(t)\! +\! B u_i(t)\! +\! L(C_1 \hat{x}_i(t)\nonumber\\
&\! -\! & \sum\limits_{j=1}^N a_{ij}(\sigma (t))(y_i(t)\! -\! y_j(t))),\ i=1,2,\ldots,N,
\end{eqnarray}
where $L\in \mathrm{R}^{n \times q_1}$ is the observer gain to be designed.

Now let $z_{tr,i}(t)\! =\! C_2(x_i(t)\! -\!\frac{1}{N}\sum\limits_{j=1}^N x_j(t)),\ i=1,2,\ldots,N,$
be the transient part of the concerned states, where $C_2 \in \mathrm{R}^{ q_2 \times n}$ is a constant matrix. With the above notation, we are ready to define the $H_{\infty}$ consensus integrated with transient performance improvement problem.

\begin{definition}
The $H_{\infty}$ consensus integrated with transient performance improvement problem is to 
design a class of controllers for the networked agents described by (\ref{dynamics_each_agent}) such that
\begin{enumerate}
\item 
when agents are not perturbed by external disturbance, i.e., $w_i(t)\! \equiv \! 0,\ i=1,2,\ldots,N$, all agents can reach consensus, i.e., $\mathbb{E}[\|x_i(t)-x_j(t)\|^2]\! = \! 0,\ i,\ j=\! 1,\ 2,\ldots,\ N$;
\item 
given $\gamma>0$, when the external disturbance exists, i.e., $w_i(t)\neq 0,\ i=1,2,\ldots,N$, the following measure index is less than $\gamma^2$, which measures the disturbance attenuation performance of the concerned output $z_{tr,i}(t),\ i=1,2,\ldots,N $:
\begin{equation}\label{measure_function}
J_{tr}\! \triangleq \!\sup\limits_{\|\omega(t)\|^2\! +\! x_0^T \bar{R} x_0\! +\! \hat{x}_0^T \bar{R} \hat{x}_0\neq 0} \frac{\mathbb{E}[z_{tr}^T(t)z_{tr}(t)]}{\|\omega(t)\|^2\! +\! x_0^T \bar{R} x_0\! +\! \hat{x}_0^T \bar{R} \hat{x}_0},
\end{equation}
where $z_{tr}(t)\! =  \! [\begin{smallmatrix}z_{tr,1}^T(t) & z_{tr,2}^T(t) & \ldots & z_{tr,N}^T(t)\end{smallmatrix}]^T$, $\omega(t)=[\begin{smallmatrix} \omega_1^T(t) &\omega_2^T(t)& \ldots & \omega_N^T(t) \end{smallmatrix}]^T$, $x_0\! = \! [\begin{smallmatrix} x_1^T(0) & x_2^T(0) & \ldots & x_N^T(0)\end{smallmatrix}]^T$ and $\hat{x}_0\! = \! [\begin{smallmatrix} \hat{x}_1^T(0) & \hat{x}_2^T(0) & \ldots & \hat{x}_N^T(0)\end{smallmatrix}]^T$; $\bar{R}\triangleq I_N \otimes R$, where $\otimes$ denotes the Kronect product and $R$ is a previously given positive-definite matrix.
\end{enumerate}
\end{definition}

Based on the observed state information, we introduce a class of controllers for each agent, which is in the following form: 
\begin{equation}\label{FS_controller_each_agent}
u_i(t)=-\tau K \hat{x}_i(t),\ i=1,2,\ldots,N,
\end{equation}
where $\tau\in \mathrm{R}$ and $K\in\mathrm{R}^{m \times n}$ are to be designed.
Hereafter, we call controller (\ref{FS_controller_each_agent}) the full-state observer-based consensus protocol, since it is constructed based on the observer (\ref{FS_observer_each_agent}).

\begin{remark}
The full-state observer-based consensus protocol (\ref{FS_controller_each_agent}) solely relies on the relative outputs of neighbouring agents; hence it is distributed.
\end{remark}

Let $x(t)\triangleq [\begin{smallmatrix}x_1(t) & x_2(t) &\ldots & x_N(t)\end{smallmatrix} ]$ and $\hat{x}(t)\triangleq [\begin{smallmatrix}\hat{x}_1(t) & \hat{x}_2(t) &\ldots & \hat{x}_N(t)\end{smallmatrix} ]$. Substituting (\ref{FS_controller_each_agent}) into (\ref{dynamics_each_agent}) and rewriting them in a compact form yield
\begin{equation}\label{FS_agent_dynamics_compact_form}
\dot{x}(t)=(I_N \otimes A)x(t)-\tau (I_N\otimes B K)\hat{x}(t) + (I_N \otimes D)\omega(t).
\end{equation}
Similarly, we have the trajectory of $\hat{x}(t)$ yielding
\begin{equation}\label{FS_agent_observer_compact_form}
\dot{\hat{x}}(t)=(I_N\otimes (A+L C_1- \tau B K))\hat{x}-(\mathcal{L}_{\sigma(t)}\otimes L C_1)x.
\end{equation}

Now, we introduce the consensus error variable ${\zeta}_i(t)\triangleq x_i(t)-\frac{1}{N}\sum\limits_{j=1}^N x_j(t)$ for each agent $i=1,2,\ldots,N$, which denotes the distance of  agent $i$'s states from the average of all agents'. Note that $\zeta(t)\triangleq [\begin{smallmatrix}\zeta_1^T(t) & \zeta_2^T(t) & \ldots & \zeta_N^T(t)\end{smallmatrix}]^T=(\mathcal{M} \otimes I_n) x(t)$, where $\mathcal{M}\triangleq I_N-\frac{1}{N} \mathbf{1}\mathbf{1}^T$. Also, let ${\delta}_i(t)\triangleq (\hat{x}_i(t)-x_i(t))-\frac{1}{N}\sum\limits_{j=1}^N (\hat{x}_j(t)-x_j(t))$ for each agent $i=1,2,\ldots,N$. We also have $\delta(t)\triangleq [\begin{smallmatrix}\delta_1^T(t) & \delta_2^T(t) & \ldots & \delta_N^T(t)\end{smallmatrix}]^T=(\mathcal{M} \otimes I_n) (\hat{x}(t)-x(t))$.

Write $\zeta(t)$ and $\delta(t)$ in a column vector and denote it by $e (t)\triangleq [\begin{smallmatrix} \zeta^T(t) & \delta^T(t)\end{smallmatrix}]^T$. With equations (\ref{FS_agent_dynamics_compact_form}) and (\ref{FS_agent_observer_compact_form}), we can write the evolution trajectory of $e(t)$ into a compact form, i.e.,
\begin{equation}\label{FS_closed_loop_dynamics}
\dot{e}(t)\! = \! F(t) e(t)+ \bar{D}\omega(t)
\end{equation}
where $F(t)=\left[\begin{matrix}I_N\otimes (A-\tau B K) & -\tau I_N\otimes (B K)\\ (I_N-\tau L_{\sigma(t)})\otimes (L C_1) & I_N \otimes (A+LC_1)\end{matrix}\right]$, and $\bar{D}=\left[\begin{matrix} I_N \otimes D\\ 0 \end{matrix}\right]$.

\subsection{Controller Design}\label{sec_main_one_two}
To address the $H_{\infty}$ consensus control integrated with transient performance improvement problem, we design the following algorithm to construct the observer (\ref{FS_observer_each_agent}) and the distributed controller (\ref{FS_controller_each_agent}):
\begin{algorithm}
\begin{description}
\item[Step 1] Given $\gamma>0$, solve the following inequalities to get a pair of feasible solution $(P_1,r_1)$
\begin{equation}\label{FS_condition_one_one}
\mbox{\scriptsize $\left[\begin{array}{*{20}{c}}
P_1A^T\! +\! AP_1\! -\!  r_1 BB^T  & D    &    P_1C_1^T       &    P_1C_2^T\\
D^T                        & -\frac{\gamma^{2}}{2} I_{\ell} & 0  &    0\\
C_1P_1                         &0        & -\frac{1}{\rho(1+\bar{\pi}^2\lambda_{\max})}I_{q_1} &   0\\
C_2P_1  &0                        &0               &   -I_{q_2}\\
\end{array}
\right]<0$},
\end{equation}
\begin{equation}\label{FS_condition_one_two}
\bigg[\begin{array}{*{20}{c}}
\frac{\gamma^2}{2\kappa}  R   & I_n  \\
I_n   & P_1 \\
\end{array}
\bigg]>0,
\end{equation}
\begin{equation}\label{FS_condition_one_three}
r_1>0,
\end{equation}
where $\lambda_{\max}$ is the largest eigenvalue of matrix ${\mathcal{L}}_{\text{un}}^T{\mathcal{L}}_{\text{un}}$; $\kappa$ denotes the largest eigenvalue of matrix $\mathcal{M}^2$.
\item[Step 2] Get $(P_2,Q,r_2)$ satisfying the following inequalities:
\begin{equation}\label{FS_condition_two_one}
\mbox{\scriptsize $\left[\begin{array}{*{20}{c}}
A^TP_2\! +\! P_2 A\! +\! C_1^TQ^T \! + \! Q C_1 \! + \! r_2 P_1^{-1} BB^T P_1^{-1}  & -P_2 D & Q^T\\
-D^T P_2 & -\frac{\gamma^{2}}{2} I_{\ell} & 0\\
Q & 0 & \frac{\rho}{8} I_n
\end{array}
\right]<0$},
\end{equation}
\begin{equation}\label{FS_condition_two_two}
P_2<\frac{\gamma^2}{2\kappa} R ,
\end{equation}
\begin{equation}\label{FS_condition_two_three}
r_2>0,
\end{equation}
where $P_1$ has been obtained from Step $1$.
\item[Step 3] Choose a small enough positive constant $\rho$ such that $\rho\! < \! 2$. If $\frac{r_1}{2\! -\! \rho}\! < \! \frac{\rho r_2}{4}$, choose the coupling strength $\tau$ to satisfy $\frac{r_1}{2\! -\! \rho}\! < \! \tau \!< \! \frac{\rho r_2}{4}$, and let $K\! = \! B^TP_1^{-1}$ and $L\! = \! P_2^{-1}Q$, otherwise go back to Step 1.
\end{description}
\end{algorithm}

Next we will show that the full-order observer-based consensus protocol  (\ref{FS_controller_each_agent}) designed according to Algorithm $1$ enables all agents to solve the  $H_{\infty}$ consensus control integrated with transient performance improvement problem.

\begin{theorem}\label{FS_theorem_one}
Under Assumptions $1$ and $2$, given a constant $\gamma>0$,
the full-order observer-based consensus protocol  (\ref{FS_controller_each_agent}) enables the group of agents described by (\ref{dynamics_each_agent}) to solve the $H_{\infty}$ consensus control integrated with transient performance improvement problem if the control gains $K,\ L$ and the coupling gain $\tau$ are chosen as in Algorithm $1$.
\end{theorem}

\textbf{Proof.}
Consider the following Lyapunov function candidate
\begin{equation}\label{lyapunov_candidate}
V(t)\! = \! \sum\limits_{i=1}^s V_i(t)
\end{equation}
where $V_i(t)\! = \! \mathbb{E}[e^T(t) \text{diag}(I_N \otimes P_1^{-1} ,I_N \otimes P_2)e(t) 1_{\sigma(t)\!=\! i}]$, in which matrices $P_1$ and $P_2$ are positive-definite.

Given $i=1,2,\ldots,s$, taking the derivative of function $V_i(t)$ along the trajectory (\ref{FS_closed_loop_dynamics}) yields 
\begin{eqnarray}\label{derivative_lyapunov_one}
\text{d}{V}_i(t)& \! = \! &  \mathbb{E}[\text{d}{e}^T(t) \text{diag}(I_N \otimes P_1^{-1} ,I_N \otimes P_2) e(t) 1_{\sigma(t)\!=\! i} ]\nonumber\\
& \! + \! & \sum\limits_{i=1}^s \mathbb{E}[{e}^T(t) \text{diag}(I_N \otimes P_1^{-1} ,I_N \otimes P_2) \text{d}{e}(t) 1_{\sigma(t)\!=\! i} ] \nonumber\\
& \! + \! &  \sum\limits_{j=1}^s q_{ji} V_j(t)\! + \! o(\text{d}t). 
 \end{eqnarray}

Since $\sigma(t)$ has a unique invariant distribution $\pi=[\begin{matrix}\pi_1 & \pi_2 & \ldots & \pi_s\end{matrix}]^T$ and $\pi_i\geq \bar{\pi}$, for all $i\in \mathcal{S}$, we have
\begin{eqnarray}\label{derivative_lyapunov_two}
\dot{V}(t) &\! \leq \! &   \mathbb{E}
[
(
\left[\begin{matrix}I_N\otimes (A-\tau B K) & -\tau I_N\otimes (B K)\\ (I_N-\tau\bar{\pi} \mathcal{L}_{\text{un}})\otimes (L C_1) & I_N \otimes (A+LC_1)\end{matrix}\right] 
\left[\begin{matrix}  \zeta(t) \\ \delta(t)
\end{matrix}\right] \nonumber\\
&\! + \! &
\left[\begin{matrix} I_N \otimes D\\ 0 \end{matrix}\right]
 \omega(t) )^T  
 \left[\begin{matrix}I_N \otimes P_1^{-1}  & 0\\ 0 & I_N \otimes P_2 \end{matrix}\right] 
  \left[\begin{matrix}  \zeta(t) \\ \delta(t)
\end{matrix}\right]  \nonumber\\
& \! + \! & 
 \left[\begin{matrix}  \zeta(t) \\ \delta(t)
\end{matrix}\right] ^T
\left[\begin{matrix}I_N \otimes P_1^{-1}  & 0\\ 0 & I_N \otimes P_2 \end{matrix}\right] \nonumber\\
& & (
\left[\begin{matrix}I_N\otimes (A-\tau B K) & -\tau I_N\otimes (B K)\\ (I_N-\tau\bar{\pi} \mathcal{L}_{\text{un}})\otimes (L C_1) & I_N \otimes (A+LC_1)\end{matrix}\right] 
\left[\begin{matrix}  \zeta(t) \\ \delta(t)
\end{matrix}\right]\nonumber\\
&\! + \!  &
\left[\begin{matrix} I_N \otimes D\\ 0 \end{matrix}\right]
 \omega(t)
)
] \nonumber \\
  & \! = \! & \mathbb{E}[\zeta^T(t)(I_N \otimes (A^TP_{1}^{-1}\! + \! P_{1}^{-1} A \! - \!  2\tau P_{1}^{-1} B B^T P_{1}^{-1}))\zeta(t)\nonumber\\
 & \!  + \! & 2  \zeta^T(t)(I_N \otimes (-\tau P_{1}^{-1} B B^T P_{1}^{-1}\! + \! C_1^T L^T P_2)\nonumber\\
 &\! + \! & \bar{\pi} \mathcal{L}_{\text{un}}^T \otimes (- C_1^T L^T P_2))\delta(t)\nonumber\\
 &\! + \! & \delta^T(t) (I_N \otimes ((A\! + \! LC_1)^T P_{2} + P_{2}(A\! + \! L C_1)))\delta(t)\nonumber\\
&\! + \! & 2 \zeta^T(t)(\mathcal{M} \otimes (P_{1}^{-1} D))\omega(t) \nonumber\\
&\! + \! & 2 \delta^T(t)(\mathcal{M} \otimes (P_{2} D))\omega(t) 
],
\end{eqnarray}
where the last equality holds by letting $K\! = \! B^T P_{1}^{-1}$.

By using Lemma $\ref{lemma_Yang}$, it gives that
\begin{eqnarray}\label{FS_yang_application_one}
&& 2  \zeta^T(t)(I_N \otimes (-\tau P_{1}^{-1} B B^T P_{1}^{-1})\delta(t)\nonumber\\
&= & 2\tau \zeta^T(t)(I_N \otimes(- P_{1}^{-1} B))(I_N \otimes(B^T P_{1}^{-1}))\delta(t)
\nonumber\\
&\leq &   \rho \tau \zeta^T(t)(I_N\otimes (P_{1}^{-1} BB^T P_{1}^{-1}))\zeta(t)\nonumber\\
& + & \frac{4}{\rho} \tau  \delta^T(t)(I_N \otimes (P_{1}^{-1} B B^T P_{1}^{-1}))\delta(t),
\end{eqnarray}
where $\rho$ is a positive constant satisfying $\rho\! < \! 2$.
Similarly, we have
\begin{eqnarray}\label{FS_yang_application_two}
&& 2  \zeta^T(t)(I_N \otimes ( C_1^T L^T P_2))\delta(t)\nonumber\\
&\leq &   \rho  \zeta^T(t)(I_N\otimes (C_{1}^{T}  C_{1}))\zeta(t)\nonumber\\
& + & \frac{4}{\rho}  \delta^T(t)(I_N \otimes (P_{2} L L^T P_{2}))\delta(t),
\end{eqnarray}
and
\begin{eqnarray}\label{FS_yang_application_three}
&& 2  \bar{\pi} \zeta^T(t)(\mathcal{L}_{\text{un}}^T \otimes (-C_1^T L^T P_2)\delta(t)\nonumber\\
&\leq &   \rho  \bar{\pi}^2 \zeta^T(t)((\mathcal{L}_{\text{un}}^T \mathcal{L}_{\text{un}})\otimes (C_{1}^{T}  C_{1}))\zeta(t)\nonumber\\
& + & \frac{4}{\rho}   \delta^T(t)(I_N \otimes (P_{2} L L^T P_{2}))\delta(t).
\end{eqnarray}

Substituting inequalities (\ref{FS_yang_application_one}), (\ref{FS_yang_application_two}) and (\ref{FS_yang_application_three}) into (\ref{derivative_lyapunov_two}) yields
\begin{eqnarray}\label{derivative_lyapunov_three}
\dot{V}(t) &\! \leq \! &
\mathbb{E}[\zeta^T(t)(I_N \otimes (A^TP_{1}^{-1}\! + \! P_{1}^{-1} A \nonumber\\
&\! - \! & \tau  (2\! - \! \rho) P_{1}^{-1} B B^T P_{1}^{-1}\!
+ \! \rho(1+\bar{\pi}^2 \lambda_{\max})C_1^T C_1))\zeta(t)\nonumber \\
&\! + \! & \delta^T(t) (I_N \otimes ((A\! + \! LC_1)^T P_{2} + P_{2}(A\! + \! L C_1)\nonumber\\
& \! + \! & \frac{4}{\rho}\tau P_{1}^{-1}  B B^T P_{1}^{-1}
\! + \! \frac{8}{\rho} P_2 L L^T P_2))\delta(t)\nonumber\\
&\! + \! & 2 \zeta^T(t)(\mathcal{M} \otimes (P_{1}^{-1} D))\omega(t) \nonumber\\
&\! + \! & 2 \delta^T(t)(\mathcal{M} \otimes (P_{2} D))\omega(t) 
].
\end{eqnarray}

Let $\tilde{\zeta}(t)\! = \! (U^T \otimes I_n)\zeta(t), \ \tilde{\delta}(t)\! = \! (U^T \otimes I_n) \delta(t)$ and $\tilde{\omega}(t) \! = \! (U^T \otimes I_n) \omega(t)$, where $U$ is a unitary matrix of matrices $\mathcal{M}$ such that $U^T \mathcal{M} U = \mathrm{diag}\{0,1,\ldots,1 \}$, where the first column of $U$ is set as $\frac{\mathbf{1}}{\sqrt{N}}$.
With the above variable changes, we have
\begin{eqnarray}\label{derivative_lyapunov_four}
\dot{V}(t) 
 &\! \leq  \! &
\sum\limits_{i=2}^N \mathbb{E}[\tilde{\zeta}_i^T(t) (A^TP_{1}^{-1}\! + \! P_{1}^{-1} A \nonumber\\
&\! - \! & \tau  (2\! - \! \rho) P_{1}^{-1} B B^T P_{1}^{-1}\!
+ \! \rho({1+\bar{\pi}^2 \lambda_{\max}})C_1^T C_1)\tilde{\zeta}_i(t)\nonumber \\
&\! + \! & \tilde{\delta}_i^T(t)  ((A\! + \! LC_1)^T P_{2} + P_{2}(A\! + \! L C_1)\nonumber\\
& \! + \! & \frac{4}{\rho}\tau P_{1}^{-1} B B^T P_{1}^{-1}
\! + \! \frac{8}{\rho} P_2 L L^T P_2)\tilde{\delta}_i(t)\nonumber\\
&\! + \! & 2 \tilde{\zeta}_i^T(t) P_{1}^{-1} D \tilde{\omega}_i(t)\! + \!  2 \tilde{\delta}_i^T(t) P_{2} D\tilde{\omega}_i(t) 
],
\end{eqnarray}
where to get the last inequality, we have used the fact that $\tilde{\zeta}_{1}(t)=\frac{\mathbf{1}^T}{\sqrt{N}}\delta(t)\equiv 0$ and
$\lambda_{\max} \triangleq \max\limits_{i=1,\ldots,N}\lambda_i(\mathcal{L}_{\text{un}}^T \mathcal{L}_{\text{un}})$.

When $\omega(t)\! \equiv \! 0$, also $\tilde{\omega}(t) \! \equiv \! 0$, we have
\begin{eqnarray}\label{derivative_lyapunov_five}
\dot{V}(t) &\! \leq  \! & 
\sum\limits_{i=2}^N \mathbb{E}[\tilde{\zeta}_i^T(t) (A^TP_{1}^{-1}\! + \! P_{1}^{-1} A \nonumber\\
&\! - \! & \tau  (2\! - \! \rho) P_{1}^{-1} B B^T P_{1}^{-1}\!
+ \! \rho(1+\bar{\pi}^2 \lambda_{max})C_1^T C_1)\tilde{\zeta}_i(t)\nonumber \\
&\! + \! & \tilde{\delta}_i^T(t)  ((A\! + \! LC_1)^T P_{2} + P_{2}(A\! + \! L C_1)\nonumber\\
& \! + \! & \frac{4}{\rho}\tau P_{1}^{-1} B B^T P_{1}^{-1}
\! + \! \frac{8}{\rho} P_2 L L^T P_2)\tilde{\delta}_i(t) ].
\end{eqnarray}
If 
\begin{eqnarray}\label{FS_deri_inequa_one_one}
&& A^TP_{1}^{-1}\! + \! P_{1}^{-1} A \! - \!  \tau  (2\! - \! \rho) P_{1}^{-1} B B^T P_{1}^{-1} \nonumber\\
&\!
+ \! &{\rho}(1+\bar{\pi}^2 \lambda_{\max})C_1^T C_1<0
\end{eqnarray}
and
\begin{eqnarray}\label{FS_deri_inequa_one_two}
&& (A\! + \! LC_1)^T P_{2} + P_{2}(A\! + \! L C_1)\nonumber\\
& \! + \! & \frac{4}{\rho}\tau P_{1}^{-1} B B^T P_{1}^{-1}
\! + \! \frac{8}{\rho} P_2 L L^T P_2<0,
\end{eqnarray}
we have $\dot{V}(t) \! \leq  \! 0$, which implies that all agents can reach consensus when the external disturbance does not exist. Next we will prove that Algorithm $1$ guarantees that inequalities (\ref{FS_deri_inequa_one_one}) and (\ref{FS_deri_inequa_one_two}) are always true. In light of Step $1$ in Algorithm $1$, we get
\begin{equation}\label{FS_deri_inequa_one_three}
P_1A^T\! + \! AP_1-r_1BB^T\! + \! {\rho}(1+\bar{\pi}^2 \lambda_{\max}) P_1 C_1^T C_1 P_1 <0.
\end{equation}
Multiplying both sides of (\ref{FS_deri_inequa_one_three}) by $P_1^{-1}$, considering that $r_1\! <\! \tau  (2 \! - \! \rho )$, and using Schur's Complement Lemma we have inequality (\ref{FS_deri_inequa_one_one}) holds. Then, substituting $Q\! = \! P_2L$ into (\ref{FS_condition_two_one}) and using the assertion that (\ref{FS_condition_one_two}) and $r_2\! > \! \frac{4}{\rho} \tau  $ yield that (\ref{FS_deri_inequa_one_two}) is true. Thus, the first condition in Definition $1$ is satisfied, which implies that the full-order observer-based control protocol enables the agents described by (\ref{dynamics_each_agent}) to reach consensus in the asymptotic mean square sense when they are not perturbed by external disturbances.

When $\omega(t) \! \neq \! 0$, we note that the last inequality of (\ref{derivative_lyapunov_four}) still holds by adding $\sum\limits_{i=2}^N \mathbb{E}[\tilde{z}_{tr,i}^T(t)\tilde{z}_{tr,i}(t) ] \! - \!\gamma^2 \sum\limits_{i=2}^N \mathbb{E}[\tilde{\omega}_i^T(t)\tilde{\omega}_i(t)]$ to the right-hand side of it and then subtracting it. In light of this, we can get the following inequality after some mathematical calculations:
\begin{eqnarray}\label{derivative_lyapunov_six}
\dot{V}(t) &\! \leq  \! & \sum\limits_{i=2}^N \mathbb{E}[\tilde{\xi}_i^T(t)\Sigma_1\tilde{\xi}_i(t)
 \! + \!  \tilde{\delta}_i^T(t)\Sigma_2 \tilde{\delta}_i(t)]\nonumber\\
& \! - \! & \sum\limits_{i=2}^N \mathbb{E}[\tilde{z}_{tr,i}^T(t)\tilde{z}_{tr,i}(t) ] \! + \! \gamma^2\sum\limits_{i=2}^N \mathbb{E}[\tilde{\omega}_i^T(t)\tilde{\omega}_i(t)], \nonumber
\end{eqnarray}
where $\tilde{z}_{tr}(t)\! =  \! [\begin{smallmatrix}\tilde{z}_{tr,1}^T(t) & \tilde{z}_{tr,2}^T(t) & \ldots & \tilde{z}_{tr,N}^T(t)\end{smallmatrix}]$, $\tilde{z}_{\text{tr}}(t)\! = \! (U^T \otimes I_n){z}_{\text{tr}}(t)$, and $\tilde{z}_{tr,1}(t)\equiv 0$; matrix $\Sigma_1$ is 
\begin{equation}\label{derivative_matrix_one_one}
\bigg[\mbox{\footnotesize ${
\begin{array}
 {*{20}{c}}
A^T P_{1}^{-1}\! + \! P_{1}^{-1}A \! - \! \tau  (2\! - \! \rho)P_1^{-1}BB^T P_1^{-1}\! + \!  {\rho}(1+\bar{\pi}^2 \lambda_{max})C_1^T C_1 \! + \! C_2^T C_2 & P_1^{-1}D  \\
D^TP_1^{-1}   & -\frac{\gamma^2}{2} I_{\ell} \\
\end{array}} $}
\bigg];
\end{equation}
matrix $\Sigma_2$ is 
\begin{equation}\label{derivative_matrix_one_one}
\bigg[\mbox{\footnotesize ${
\begin{array}
 {*{20}{c}}
(A\! + \! LC_1)^T P_{2} + P_{2}(A\! + \! L C_1) \! + \!  \frac{4}{\rho}\tau P_{1}^{-1} B B^T P_{1}^{-1}
\! + \! \frac{8}{\rho} P_2 L L^T P_2 & -P_2D\\
-D^T P_2 & -\frac{\gamma^2}{2} I_{\ell} \\
\end{array}} $}
\bigg].
\end{equation}

If the following inequalities hold,
\begin{equation}\label{FS_deri_inequa_two_one}
\Sigma_1\! < \! 0,
\end{equation}
\begin{equation}\label{FS_deri_inequa_two_two}
\Sigma_2\! < \! 0,
\end{equation}
we have 
\begin{eqnarray}\label{derivative_lyapunov_seven}
\dot{V}(t) &\! \leq  \! & \! - \! \sum\limits_{i=2}^N \mathbb{E}[\tilde{z}_{tr,i}^T(t)\tilde{z}_{tr,i}(t) ] \! + \! \sum\limits_{i=2}^N \mathbb{E}[\tilde{\omega}_i^T(t)\tilde{\omega}_i(t)].
\end{eqnarray}
Taking similar steps to prove inequalities (\ref{FS_deri_inequa_one_one}) and (\ref{FS_deri_inequa_one_two}), we can verify that inequalities (\ref{FS_deri_inequa_two_one}) and (\ref{FS_deri_inequa_two_two}) are satisfied under Algorithm $1$.

Integrating both sides of inequality (\ref{derivative_lyapunov_seven}) from zero to infinity yields
\begin{eqnarray}\label{FS_deri_inequa_three_one}
\sum\limits_{i=2}^N \mathbb{E}[\|\tilde{z}_i(t)\|^2]& \! \leq \! &  \zeta^T(0) (I_N\otimes P_1^{-1} )\zeta(0)
 \! + \!  \delta^T(0) (I_N\otimes P_2 )\delta(0)\nonumber \\
& \! + \! & \gamma^2 \sum\limits_{i=2}^N \mathbb{E}[\tilde{\omega}_i^T(t)\tilde{\omega}_i(t)]\nonumber\\
& \! = \! & 
\zeta^T(0) (\mathcal{M}^2 \otimes P_1^{-1} )\zeta(0)
 \! + \!  \delta^T(0) (\mathcal{M}^2 \otimes P_2 )\delta(0)\nonumber \\
& \! + \! & \gamma^2 \sum\limits_{i=2}^N \mathbb{E}[\tilde{\omega}_i^T(t)\tilde{\omega}_i(t)]\nonumber\\
& \! \leq \! & 
\zeta^T(0) (\mathcal{M}^2 \otimes P_1^{-1} )\zeta(0)
 \! + \!  \delta^T(0) (\mathcal{M}^2 \otimes P_2 )\delta(0)\nonumber \\
& \! + \! & \gamma^2 \sum\limits_{i=1}^N \mathbb{E}[\tilde{\omega}_i^T(t)\tilde{\omega}_i(t)]\nonumber\\
& \! = \! & 
\zeta^T(0) (\mathcal{M}^2 \otimes P_1^{-1} )\zeta(0)
 \! + \!  \delta^T(0) (\mathcal{M}^2 \otimes P_2 )\delta(0)\nonumber \\
& \! + \! & \gamma^2\|\tilde{\omega}(t)\|^2\nonumber\\
&\leq \! & \gamma^2 x^T(0) (I_N \otimes R)x(0)
 \! + \! \gamma^2 \hat{x}^T(0) (I_N \otimes R)\hat{x}(0)\nonumber\\
& \! + \! & \gamma^2 \|\tilde{\omega}(t)\|^2,
\end{eqnarray}
where the last equality is true because the external disturbance considered herein is deterministic and thus $\sum\limits_{i=1}^N \mathbb{E}[\tilde{\omega}_i^T(t)\tilde{\omega}_i(t)]= \|\tilde{\omega}(t)\|^2$; the last inequality also holds since 
\begin{equation}\label{FS_deri_inequa_three_two}
\kappa P_1^{-1}\! < \! \frac{\gamma^2}{2}  R,
\end{equation}
and
\begin{equation}\label{FS_deri_inequa_three_three}
\kappa P_2\! < \!  \frac{\gamma^2}{2} R \! < \! \gamma^2  R.
\end{equation}

Since $\tilde{z}_{tr,1}(t)\! \equiv\! 0$, we have $\sum\limits_{i=2}^N \mathbb{E}[\|\tilde{z}_{tr,i}(t)\|^2]\! = \! \sum\limits_{i=1}^N \mathbb{E}[\|\tilde{z}_{tr,i}(t)\|^2]\! = \!  \mathbb{E}[\|{z}_{tr}(t)\|^2]$. It follows from (\ref{FS_deri_inequa_three_one}) that 
\begin{equation}\label{FS_objective_hold}
J_{\text{tr}}\! < \! \gamma^2,
\end{equation}
which means that the second condition in Definition $1$ is satisfied. That is, the $H_{\infty}$ consensus integrated with transient performance improvement problem is solved. The proof is completed.

\begin{remark}
The main advantages of the full-order observer-based consensus protocol (\ref{FS_controller_each_agent}) are two-fold. Firstly, as mentioned in Remark $1$, it solely require the relative output of neighbours which are measurable by a set of sensor networks in practice. Note that consensus protocols in \cite{Wang2015distributed,zhao2012hinfinity}, which deal with the $H_{\infty}$ consensus for general linear multi-agent systems over fixed communication graph, need the exchange of the protocol's internal states or virtual outputs, apart from the relative outputs, to maintain ``the separation principle". In fact, the states of the observers are internal information for the agents, which have to be transmitted via some communication networks. Thus, compared with protocols in \cite{Wang2015distributed,zhao2012hinfinity}, the control law (\ref{FS_controller_each_agent}) can reduce or even remove the communication burden imposed on the communication networks. Secondly, to the best of our knowledge, most of the references \cite{Mu2014L,you2013consensus,Wang2017distributed,wang2018h} related to the consensus control of multi-agent systems over Markovian randomly switching graphs assume that each possible graph is balanced. If we relax this assumption and let each possible graph be strongly connected, then there is no guarantee that protocols in \cite{Mu2014L,you2013consensus,Wang2017distributed,wang2018h} can serve their purpose. However, 
the control law (\ref{FS_controller_each_agent}) can be used to deal with this problem by slightly modifying the consensus error variable. 
\end{remark}

\subsection{Reduced-Order Observer-Based Controller}\label{sec_main_one_three}
Though the full-order observer-based consensus protocol (\ref{FS_controller_each_agent}) is effective to solve the $H_{\infty}$ consensus integrated with transient performance problem, it still possesses a certain degree of redundancy, which is caused by the fact that the entire state is estimated by the observer (\ref{FS_observer_each_agent}), but actually part of the state information is available in the output. In this subsection, we will eliminate the redundancy and thus reduce the dimension of the protocol, especially for the case that agents are described by multi-input multi-output dynamics.
To solve the problem mentioned in Section \ref{sec_main_one_one}, we introduce a reduced-order observer-based consensus protocol which is based on the relative information between neighbouring nodes. The reduced-order observer-based consensus protocol is given by
\begin{eqnarray}\label{RO_each_agent_controller}
\dot{v}_i(t)& \! = \! &\bar{F} v_i(t)\! + \! G y_i(t) \! + \! TB u_i(t) \! + \! TD \omega_i(t),\nonumber\\
 u_i(t)&\! =  \! & -\tau K[\sum\limits_{j=1}^N a_{ij}(\sigma (t))(Q_1(y_i(t)\! - \! y_j(t)) \! + \! Q_2 (v_i(t)\! - \! v_j(t)))],\nonumber\\
  i &\!  =  \! & 1,2,\ldots,N,
\end{eqnarray}
where $\bar{F}$ is a Hurwitz matrix and has no eigenvalues in common with those of matrix $A$. Matrices $G$ and $T$ are the unique solution to the following Sylvester equation
\begin{equation}\label{sylverster_equation}
TA\! - \! \bar{F}T\! = \! GC_1
\end{equation} 
which further satisfies that matrix $[\begin{smallmatrix}C_1^T & T^T \end{smallmatrix}]^T$ is non-singular. $Q_1$ and $Q_2$ are given by $[\begin{smallmatrix}Q_1 & Q_2 \end{smallmatrix}]\! = \! [\begin{smallmatrix}C_1^T & T^T \end{smallmatrix}]^T$.

Now we can present the algorithm to design the reduced-order observer-based consensus protocol (\ref{RO_each_agent_controller}).

\begin{algorithm}
\begin{description}
\item[Step 1] Choose matrices $\bar{F},\ G,\ T,\ Q_1$ and $Q_2$ to satisfy the requirements given above.
\item[Step 2] Given a positive constant $\gamma$, solve the following inequalities to get a pair of feasible solution $(P_1,r_1)$
\begin{equation}\label{RO_condition_one_one}
\bigg[\begin{array}{*{20}{c}}
P_1A^T\! +\! AP_1\! -\!  r_1 BB^T  & D              &    P_1C_2^T\\
D^T                        & -\gamma^{2} I_{\ell}  &    0\\
C_2P_1                         &0               &   -I_q\\
\end{array}
\bigg]<0,
\end{equation}
\begin{equation}\label{RO_condition_one_two}
\bigg[\begin{array}{*{20}{c}}
\frac{\gamma^2}{\kappa} R   & I_n  \\
I_n   & P_1 \\
\end{array}
\bigg]>0,
\end{equation}
\begin{equation}\label{RO_condition_one_three}
r_1>0.
\end{equation}
\item[Step 3] Get $(P_2,r_2)$ from the following inequalities:
\begin{equation}\label{RO_condition_two_one}
\bar{F}^TP_2\! +\! P_2 \bar{F}\! +\! r_2 Q_2 P_1^{-1} BB^T P_1^{-1} Q_2^T  <0,
\end{equation}
\begin{equation}\label{RO_condition_two_two}
P_2<\frac{\gamma^2}{\kappa}R,
\end{equation}
\begin{equation}\label{RO_condition_two_three}
r_2>0,
\end{equation}
where $P_1$ has been obtained from Step $1$.
\item[Step 4] Choose a small enough positive constant $\rho$ such that $\rho\! < \! \lambda_{\text{min}}$, where $\lambda_{\text{min}}$ is the second smallest eigenvalue of matrix $\tilde{\mathcal{L}}_{\text{un}}\! = \! {\mathcal{L}}_{\text{un}}^T \! + \! {\mathcal{L}}_{\text{un}}$. If $\frac{r_1}{\bar{\pi}(\lambda_{\text{min}}\! -\! \rho)}\! < \! \frac{\rho r_2}{4\bar{\pi} \lambda_{\text{max}}}$, let $\frac{r_1}{\bar{\pi}(\lambda_{\text{min}}\! -\! \rho)}\! < \! \tau \!< \! \frac{\rho r_2}{4\bar{\pi} \lambda_{\text{max}}}$, and $K\! = \! B^TP_1^{-1}$, where $\lambda_{\text{max}}$ is the largest eigenvalue of matrix ${\mathcal{L}}_{\text{un}}^T{\mathcal{L}}_{\text{un}}$; otherwise go back to Step 1.
\end{description}
\end{algorithm}

\begin{remark}
A necessary and sufficient condition for the feasibility of Algorithms $1$ and $2$ is that $(A,\ B,\ C)$ is stabilizable and detectable.

\end{remark}

\begin{theorem}\label{RO_theorem_one}
Under Assumptions $1$ and $2$, given a constant $\gamma>0$, the reduced-order observer-based consensus controller (\ref{RO_each_agent_controller}) enables agents described by (\ref{dynamics_each_agent}) to solve the $H_{\infty}$ consensus integrated with transient performance improvement problem if the control gain $K$ and the coupling gain $\tau$ are chosen as in Algorithm $2$.
\end{theorem}
\textbf{Proof.}
Theorem \ref{RO_theorem_one} can be proved by taking similar steps as in the proof of Theorem \ref{FS_theorem_one}.

\section{Simulation Results}\label{sec_simulation}
In this section, we illustrate the effectiveness of Algorithms $1$ and $2$ with the following example.

 Consider a network of four Raptor-$90$ helicopters whose dynamics are given by (\ref{dynamics_each_agent}) where $A$ is given by (\ref{simulation_agent_gain_A}) and
\begin{figure*}[ht]
\begin{equation}\label{simulation_agent_gain_A}
A=\mbox{\scriptsize $\left( {\begin{array}{*{13}{c}}
  {-0.1778} & 0 &  0&  0 & 0 & {-9.7807} &  {-9.7807} & 0 & 0 & 0 & 0\\
  0& -0.3104 &0& 0& 9.7807 &0 &0 &9.7807& 0 &0 &0\\
  {-0.3326} & {-0.5353} & 0 & 0  &0 & 0 & 75.7640 & 343.8600 & 0 & 0 & 0\\
   0.1903 & {-0.2940} & 0  & 0 & 0 & 0 & 172.6200 & {-59.9580} & 0 & 0 & 0 \\
   0 & 0 &  1 & 0 & 0 & 0 & 0 & 0 & 0 & 0 & 0\\
   0  & 0 & 0 & 1 & 0 &  0  & 0 & 0 & 0  &0  &0\\
   0 & 0 & 0 & -1 & 0 & 0 & -8.1222  & 4.6535 & 0 & 0 & 0\\
   0 & 0  & -1&  0&  0& 0& -0.0921& -8.1222 & 0 & 0 & 0\\
   0 & 0 & 0 &0 & 0 & 0& 17.1680 & 7.1018 & -0.6821 & -0.1070 & 0\\
   0 & 0 & -0.2834&  0 &0 & 0 &0 &0 &-0.1446& -5.5561 & -36.6740\\
   0& 0& 0& 0& 0 &0 &0 &0 &0 &2.7492& -11.1120
\end{array}} \right)$}
\end{equation}
\end{figure*}

\[
\setlength{\arraycolsep}{2pt}
\renewcommand{\arraystretch}{0.6}
B=\mbox{\small $\left( {\begin{array}{*{20}{c}}
   0 & 0 &0\\
   0& 0& 0\\
   0 &0 &0\\
   0 &0 &0\\
   0 &0 &0\\
   0 &0 &0\\
   0.0632& 3.3390& 0\\
   3.1739 &0.2216& 0\\
   0 &0 &0\\
   0 &0& -74.364\\
   0& 0& 0
\end{array}} \right)$},
D=\mbox{\small $\left( {\begin{array}{*{20}{c}}
   {-0.1778} & 0\\
   0& -0.3104\\
   {-0.3326} & {-0.5353}\\
   0.1903 & {-0.2940}\\
   0 & 0\\
   0 & 0\\
   0 & 0\\
   0 & 0\\
   0 & 0\\
   0 & 0\\
   0 & 0\\
\end{array}} \right)$},
\]
\[
\setlength{\arraycolsep}{3pt}
\renewcommand{\arraystretch}{0.6}
C_1=C_2=\mbox{\small $\left( {\begin{array}{*{20}{c}}
   0 & 0& 0 & 0&1 & 0 & 0& 0 & 0& 0 & 0\\
   0 & 0& 0 & 0&0 & 1  & 0& 0 & 0 & 0 & 0\\
   0 & 0& 0 & 0&0 & 0& 0 & 0 &0&  1& 0 \\
    0& 0 &  1& 0 & 0 & 0& 0 & 0 & 0& 0 & 0\\
    0 & 0 & 0& 1&0 & 0& 0 & 0 & 0& 0 & 0 \\
\end{array}} \right)$}.
\]
This linear state-space model at a hovering point is  established by Chen in \cite{report_TYAA_flying_vehicle_2005}.
 The state vector is
$\setlength{\arraycolsep}{2pt}
\renewcommand{\arraystretch}{1}
x=\mbox{\scriptsize$\left[ {\begin{array}{*{20}{c}}
   U^T & V^T & p^T & q^T & \phi^T & \theta^T & a^T_s & b^T_s  &     W^T &     r^T &     r^T_{fb}
\end{array}} \right]$}^T.
$
All the above variables are described in Table~\ref{variables}.

\begin{table}[!t]
\caption{Helicopter variables.}
\label{variables}
\begin{tabular}{c||c }
\hline
Variables &  \\
\hline
q & Pitch rate in the body frame components\\
p & Roll rate in the body frame components\\
r & Yaw rate in the body frame components\\
U & Velocity along the body frame x-axis \\
V & Velocity along the body frame y-axis \\
 W & Velocity along the body frame z-axis\\
  $\theta$ & Pitch angle\\
 $\phi$ & Roll angle\\
 $a_s$ & Longitudinal blade angle\\
  $r_{fb}$ & Yaw rate feedback\\
 $b_s$ & Lateral blade angle\\
\hline
\end{tabular}
\end{table}

\begin{figure}\label{communication_graph}
\subfloat[$\mathcal{G}_1$]{\includegraphics[width=1.5in]{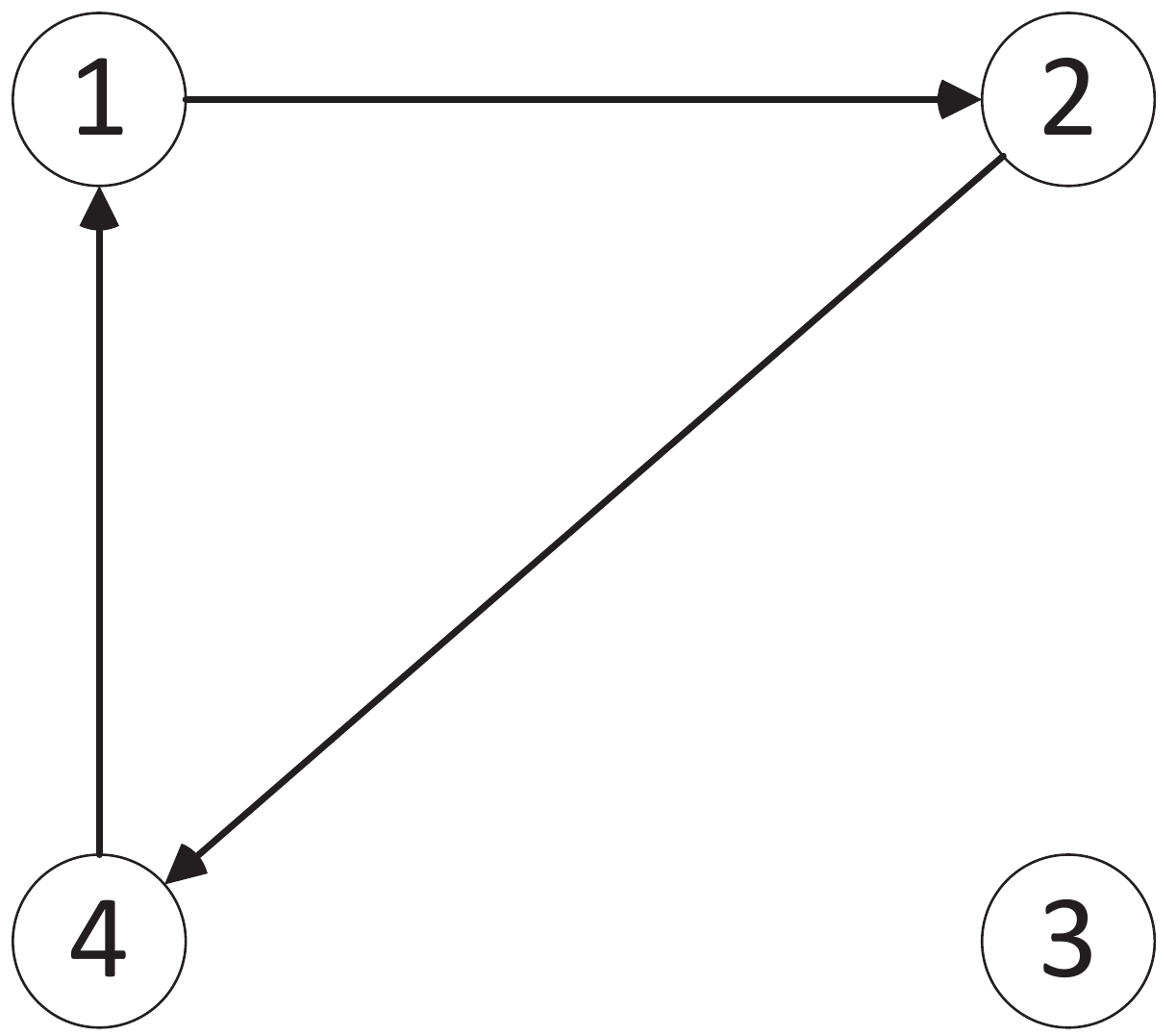}%
\label{fig_first_case}}
\hfil
\subfloat[$\mathcal{G}_2$]{\includegraphics[width=1.5in]{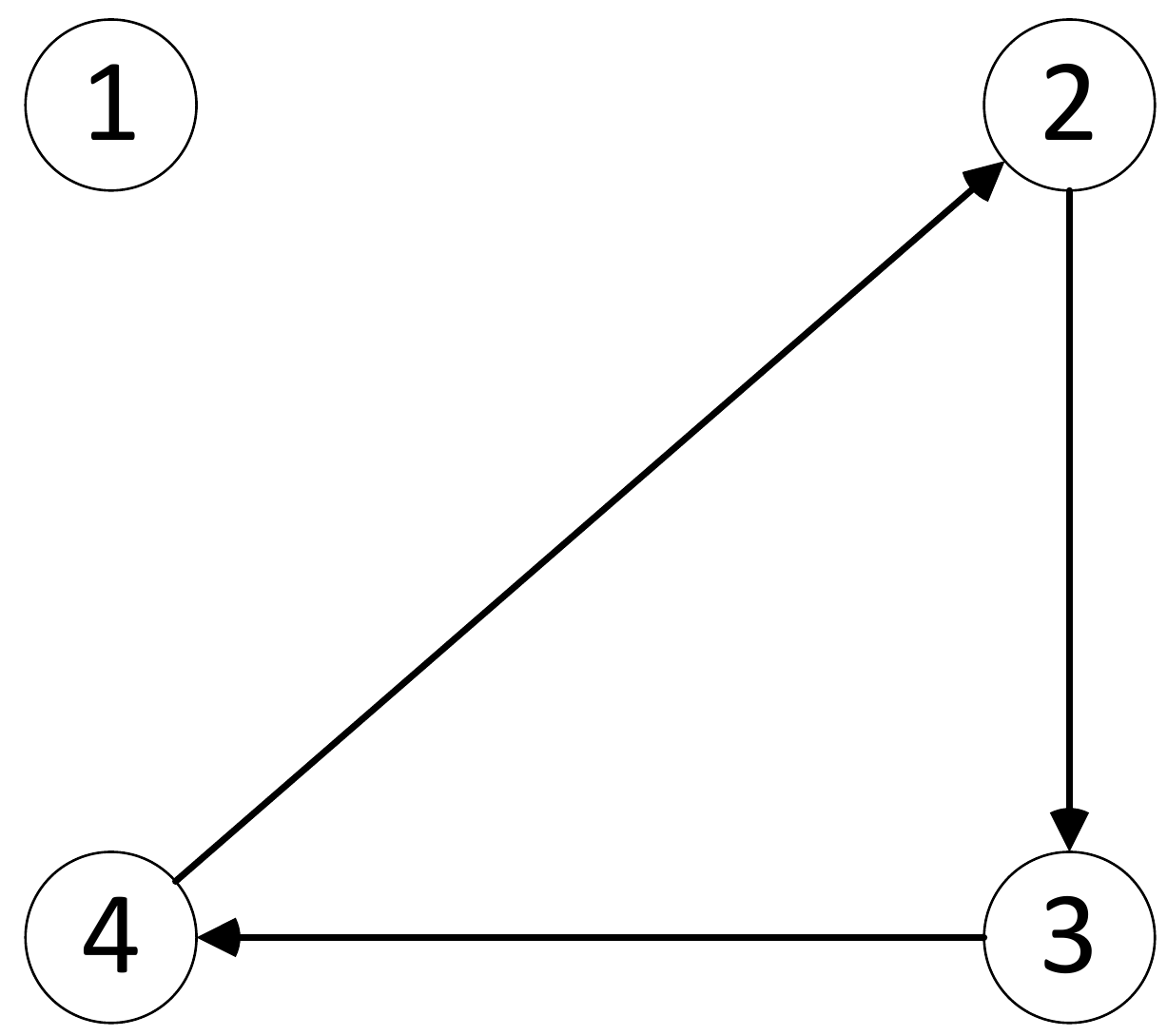}%
\label{fig_second_case}}
\hfil
\centering
\subfloat[$\mathcal{G}_{\text{un}}$]{\includegraphics[width=1.5in]{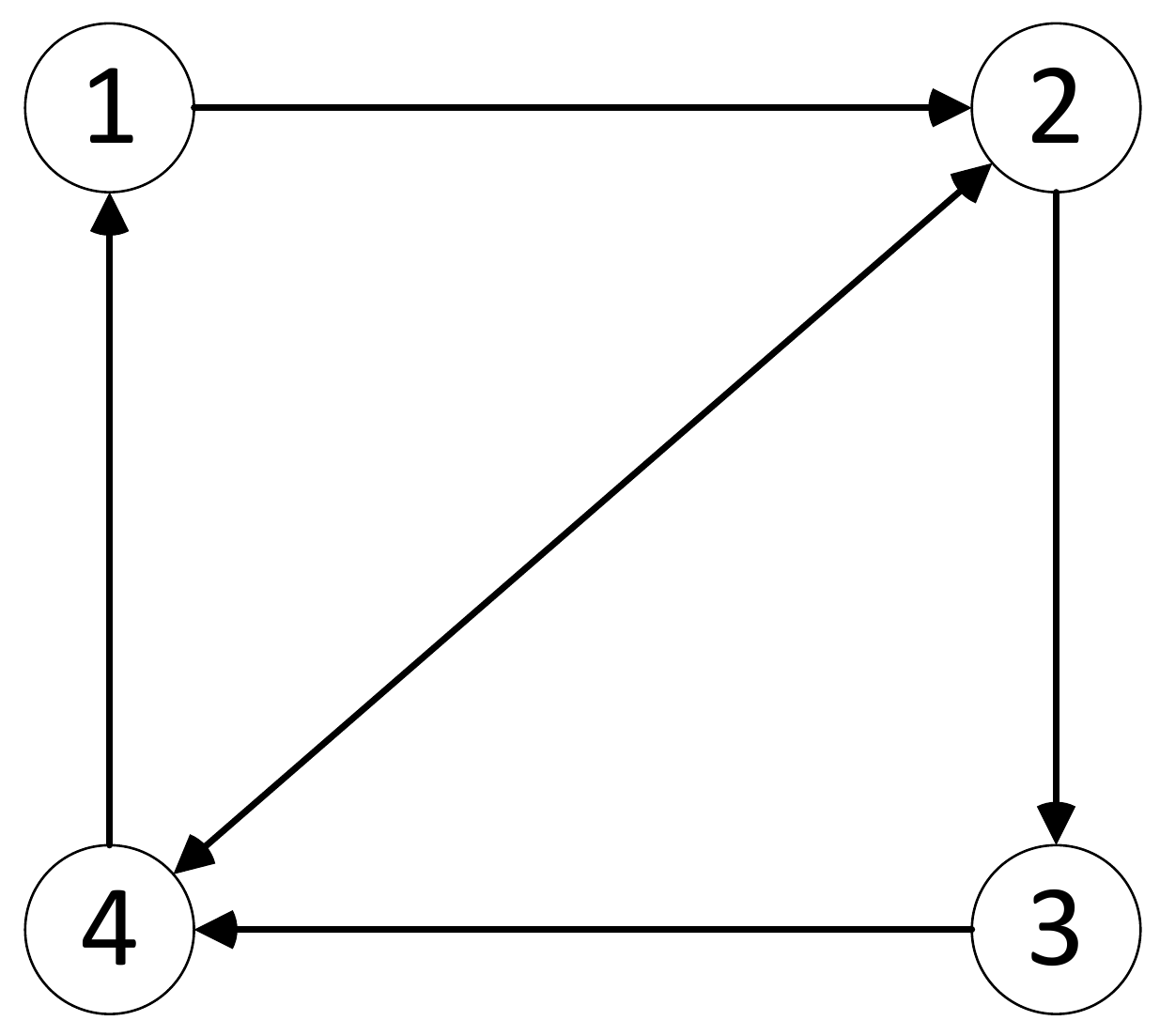}%
\label{fig_third_case}}
\caption{The communication graphs: (a) $\mathcal{G}_1$; (b) $\mathcal{G}_2$; (c) $\mathcal{G}_{\text{un}}\! = \! \mathcal{G}_1 \cup \mathcal{G}_2$.}
\end{figure}

The communication topology randomly switches between graphs $\mathcal{G}_1$ and $\mathcal{G}_2$ shown in Figs.~\ref{fig_first_case} and \ref{fig_second_case}, respectively, both of which are directed and balanced. Fig.~\ref{fig_third_case} shows the union of $\mathcal{G}_1$ and $\mathcal{G}_2$. Thus, Assumption $2$ is satisfied since the graph $\mathcal{G}_{\text{un}}$ contains a directed spanning tree. Moreover,
let the generator matrix of the continuous-time Markov process be
\begin{equation}\label{simulation_generator_matrix}
Q\! = \!\left[ \begin{array}{* {20} {c}}
-1 & 1\\
2 & -2
\end{array}\right].
\end{equation}
We assume that
the Markov process satisfies Assumption $1$ and has a unique invariant distribution $\pi\! = \! [\begin{array}{* {20} {c}}
2/3 &1/3
\end{array}  ]$.

\begin{figure*}[ht]
\begin{equation}\label{FS_infinity_transient_control_gain}
K_l=\mbox{\scriptsize $\left( {\begin{array}{*{11}{c}}
0.5701 & 2.6510 & 2.4895 & -0.9451 & 16.5745 & -3.3893 & 3.5827 & 59.2827 & 0.8457 & -0.2882 & -0.4985\\
-3.3829 & 1.0260 & 0.9998 & 3.3994 & 4.2455 & 20.5977 & 48.2969 & 6.9465 & -2.2794 & 0.3203 &1.2511\\
-14.7402 & -2.8399 & 0.0330 & -0.9026 & -14.5324 & 4.7567 & -7.5918 & 6.9036 & 0.7930 &-17.3624 &-16.0586
\end{array}} \right)$}
\end{equation}
\end{figure*}

\begin{equation}\label{FS_infinity_transient_observer_gain}
L=\mbox{\scriptsize $\left( {\begin{array}{*{5}{c}}
   -1.9899 & -4.3140 & -0.3230 &-3.4616  & 10.0917\\
     -15.0381 & 3.0932 & 4.1833 & -116.8244  &-6.1682\\
      -6.1358 & -347.6558 & 6.5180 & 35.7853 & -297.0615\\
  -10.0725 & 2.4948 & 5.9756 & 496.4880 & -7.8850\\
  -1.9442  & 3.7207 & 1.7848 & 18.6331 & -1.5905\\
  -4.0505 &-3.9455 & 5.7071 & -113.7605  & 0.0654\\
   5.3483 &-0.1549  & -5.8315 &-38.6820  &-43.9203\\
   2.1099 & -3.0888 & -5.7404 & -99.2847 & 16.7910\\
   -2.6503 &   -5.1248 & 1.6962 & 7.1855 & -6.5478\\
   9.4013 & 30.0619 & -2.2568 & -21.8266 & -12.6670\\
   2.3712 & 7.4457 & 0.6638 & -120.5388 & -3.0863
\end{array}} \right)$}
\end{equation}

\begin{equation}\label{RO_infinity_transient_control_gain_barF}
\bar{F}=\mbox{\tiny $\left( {\begin{array}{*{11}{c}}
0 & 1 & 0 & 0 & 0 & 0 \\
0 & 0 & 1 & 0 & 0 & 0\\
0 & 0 & 0 & 1 & 0 & 0\\
0 & 0 & 0 & 0 & 1 & 0\\
0 & 0 & 0 & 0 & 0 & 1\\
-8.2944 & -40.0896 & -77.2416 & -75.3600 & -38.9200 & -10
\end{array}} \right)$}
\end{equation}

\begin{equation}\label{RO_infinity_transient_control_gain_G}
G=\mbox{\scriptsize $\left( {\begin{array}{*{11}{c}}
0.3435 & 0.2485 & 0.6139 & 0.1746 & 0.3182\\
0.6631 & 0.9087 & 0.6521 & 0.0599 & 0.9556\\
0.5162 & 0.8895 & 0.6013 & 0.1524 & 0.0290\\
0.7967 & 0.9898 & 0.7978 & 0.3834 & 0.3972\\
0.5766 & 0.3237 & 0.6104 & 0.0131 & 0.2728\\
0.0669 & 0.9874 & 0.3772 & 0.9654 & 0.3619
\end{array}} \right)$}
\end{equation}

\begin{figure*}[ht]
\begin{equation}\label{RO_infinity_transient_control_gain}
K_r=\mbox{\scriptsize $\left( {\begin{array}{*{11}{c}}
   0.5701 & 2.6510 & 2.4895 &-0.9451  & 16.5745 & -3.3893 & 3.5827 & 59.2827 & 0.8457  &-0.2882  & -0.4985\\
  -3.3829 & 1.0260 & 0.9998 & 3.3994 & 4.2455  & 20.5977 & 48.2969 & 6.9465 & -2.2794 &0.3203  & 1.2511\\
   -14.7402 &-2.8399  & 0.0330 &-0.9026  &-14.5324 & 4.7567 & -7.5918 & 6.9038  & 0.7930 &   -17.3624 & -16.0856
\end{array}} \right)$}
\end{equation}
\end{figure*}

\begin{figure*}[ht]
\begin{equation}\label{RO_infinity_transient_control_gain_T}
T=\mbox{\tiny $\left( {\begin{array}{*{11}{c}}
79.3 & -314.7 & -103.7 & -351.4 & 3.2 \times {10}^3 & 1.7 \times {10}^3 & 1.4\times {10}^3 & 3.1 \times {10}^3 & 4.7 \times {10}^3 &-34.8 & 122.2\\
-47.4 & 256.3 & 48.9 & 268.6 & -3.0 \times {10}^3 & -776.5 & -587.7 & -3.0 \times {10}^3 & -3.2 \times {10}^3 & 23.6 & -83.2  \\
41.9 & -185.6 & -23.9 & -189.7 & 2.5 \times {10}^3 & 463.1 & 337.8 & 2.4 \times {10}^3 & 2.1\times {10}^3 &-16.1 & 56.7\\
-36.7 & 125.3 & 14.4 & 125.3 & -1.8\times {10}^3 & -411.4 & -330.3 & -1.8\times {10}^3 & -1.4\times {10}^3 & 10.9 & -38.5\\
23.9 & -84.4 & -10.4 & -81.4 & 1.2\times {10}^3 & 358.2 & 301.7 & 1.2\times {10}^3 & 1.0\times {10}^3 & -7.6 & 26.2\\
-17.4 & 55.4  & 7.5 & 56.2 &-826.5 & -235.0 & -200.2 & -822.2 & -692.5 & 5.2 & -11.3
\end{array}} \right)$}
\end{equation}
\end{figure*}

\begin{equation}\label{RO_infinity_transient_control_gain_Q1}
Q_1=\mbox{\tiny $\left( {\begin{array}{*{11}{c}}
-0.0241 & 0.015   & 1.0730  & -0.0456 & -0.1051\\
-0.0037 & 0.6453  & -0.4305 & 0.0460  & 0.1019\\
0.0368  & 0.4588  & 1.1349  & 0.9725  & -0.0056\\
0.2473  & -0.5468 & -0.6403 & 0.0862  & 1.1512\\
0.2296  & -0.6722 & -1.4509 & 0.0494  & 0.1288\\
0.0247  & 0.4824  &1.1776   & -0.0276 & -0.0637\\
-0.0368 & -0.4588 & -1.1349 & 0.0275  & 0.0056\\
-0.2264 & 0.6745  & 1.3874  & -0.0678 & -0.1130\\
0.0241  & -0.0105 & -0.0730 & 0.0456  & 0.1051\\
0.7944  & 0.6617  & 0.3779  & -0.0037 & -0.0237\\
-0.0209 & -0.1277 & -0.7471 & -0.0184 & -0.0382
\end{array}} \right)$}
\end{equation}

\begin{equation}\label{RO_infinity_transient_control_gain_Q2}
Q_2=\mbox{\tiny $\left( {\begin{array}{*{11}{c}}
-0.0085 & -0.0389 & -0.0747 & -0.0716 & -0.0144 & 0.0209\\
0.0049  & 0.0163  & 0.0065  & -0.0348 & -0.1162 & -0.1174\\
-0.0181 & -0.0465 & 0.0168  & 0.2096  & 0.1630  & -0.0654\\
0.1151  & 0.5554  & 1.0613  & 1.0445  & 0.3863  & -0.1240\\
0.1371  & 0.6288  & 1.1174  & 0.9338  & 0.3509  & 0.0482\\
-0.0250 & -0.0795 & -0.0464 & 0.1556  & 0.1525  & -0.0582\\
0.0181  & 0.0465  &-0.0168  & -0.2096 & -0.1630 & 0.0654\\
-0.1351 & -0.6186 & -1.1013 & -0.9237 & -0.3500 & -0.0516\\
0.0085  & 0.0389  & 0.0747  & 0.0716  & 0.0144  & -0.0209\\
-0.1286 & -0.5899 & -1.0427 & -0.8622 & -0.3365 & -0.0691\\  
0.0200  & 0.0632  & 0.0400  & -01209  & -0.0363 & 0.1757
\end{array}} \right)$}
\end{equation}

Choose $R=I$ and
let $\gamma=4$. Using the Yalmip toolbox, we can obtain the control gain matrices $K_l$ and $L$ of the full-state observer-based consensus protocol through Algorithm $1$, given by (\ref{FS_infinity_transient_control_gain}) and (\ref{FS_infinity_transient_observer_gain}), respectively. Since Algorithm $1$ has the feature of decoupling the design of the agent dynamics from the communication topologies, we can select the coupling gain by following Step $3$ of Algorithm $1$.
And we set it as $1.4$ in this example.

Also, choosing the same $R$ and $\gamma$, we can get the control protocol gain matrices $\bar{F}$, $G$, $K_r$, $T$, $Q_1$, and $Q_2$ of the reduced-order observer-based consensus protocol through Algorithm $2$, which are presented in (\ref{RO_infinity_transient_control_gain_barF}),
(\ref{RO_infinity_transient_control_gain_G}), (\ref{RO_infinity_transient_control_gain}), (\ref{RO_infinity_transient_control_gain_T}), (\ref{RO_infinity_transient_control_gain_Q1}), and (\ref{RO_infinity_transient_control_gain_Q2}), respectively. The coupling gain is also taken as $1.4$.

\begin{figure}
\includegraphics[width=4in]{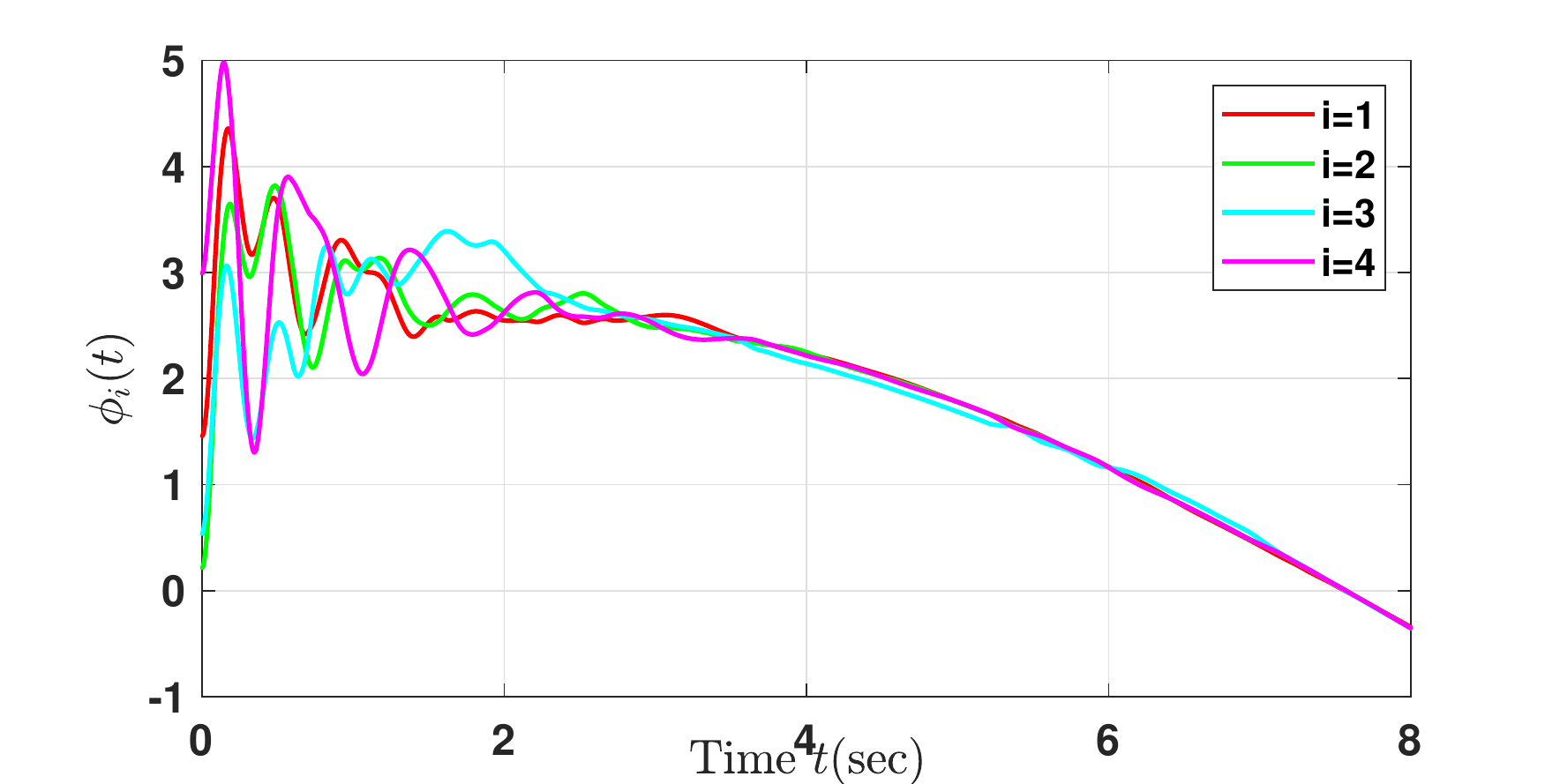}
\caption{The roll angles of four agents under the full-state observer-based consensus protocol (\ref{FS_controller_each_agent}).}
\label{fig_phi_fs_trans}
\end{figure}

\begin{figure}
\includegraphics[width=4in]{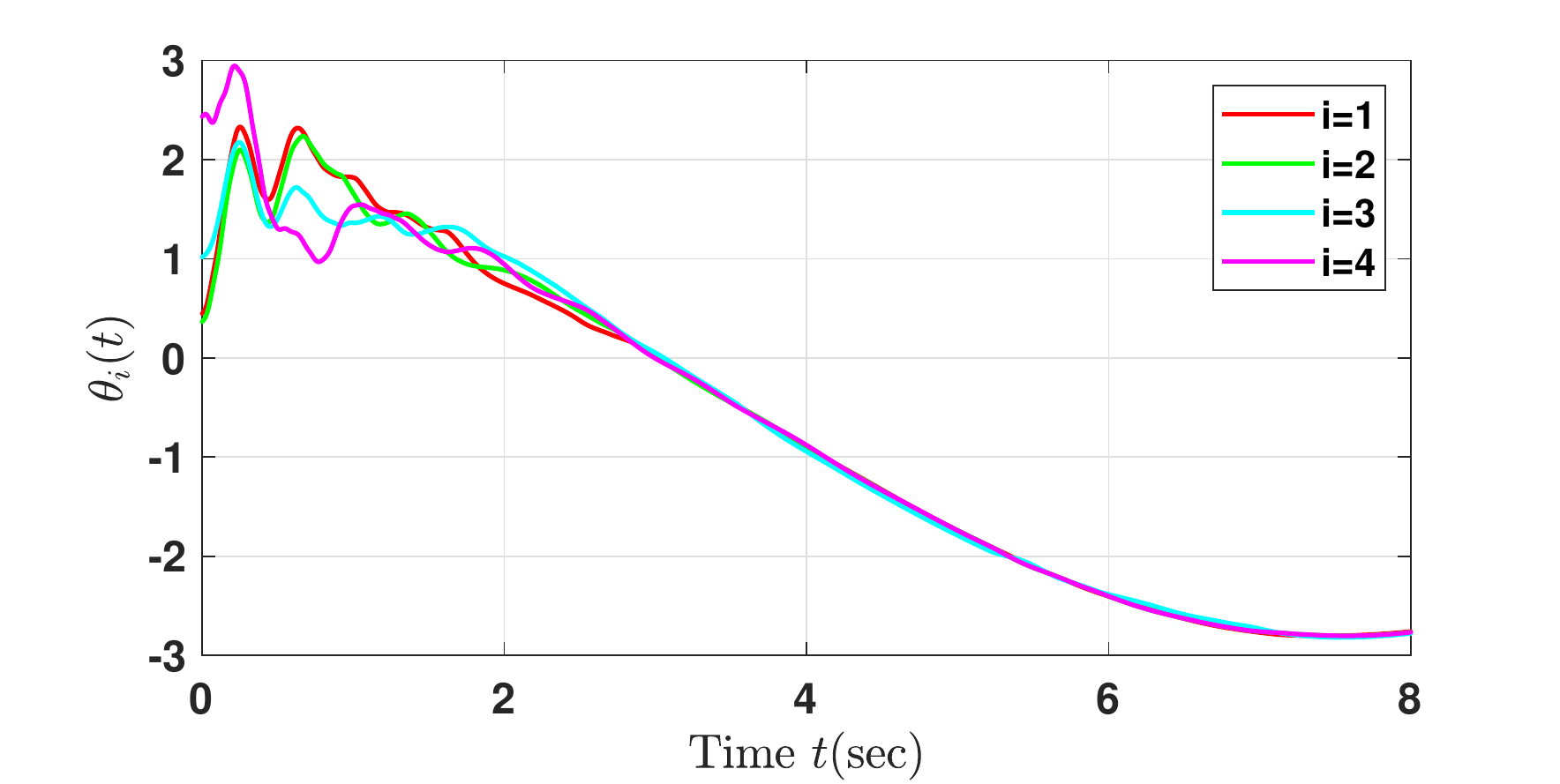}
\caption{The pitch angles of four agents under the full-state observer-based consensus protocol (\ref{FS_controller_each_agent}).}
\label{fig_theta_fs_trans}
\end{figure}

\begin{figure}
\includegraphics[width=4in]{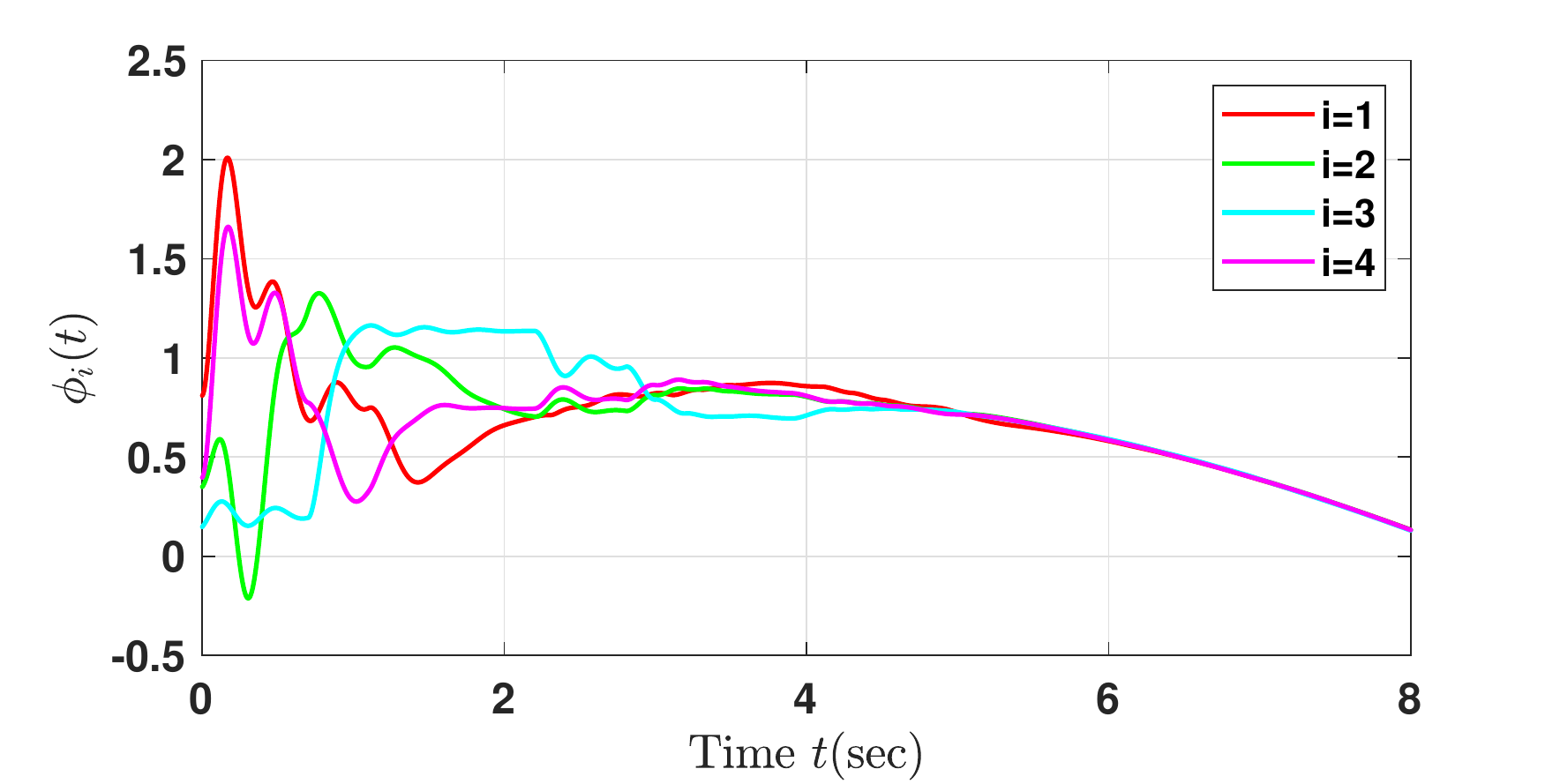}
\caption{The roll angles of four agents under the reduced-order observer-based consensus protocol (\ref{RO_each_agent_controller}).}
\label{fig_phi_ro_trans}
\end{figure}

\begin{figure}
\includegraphics[width=4in]{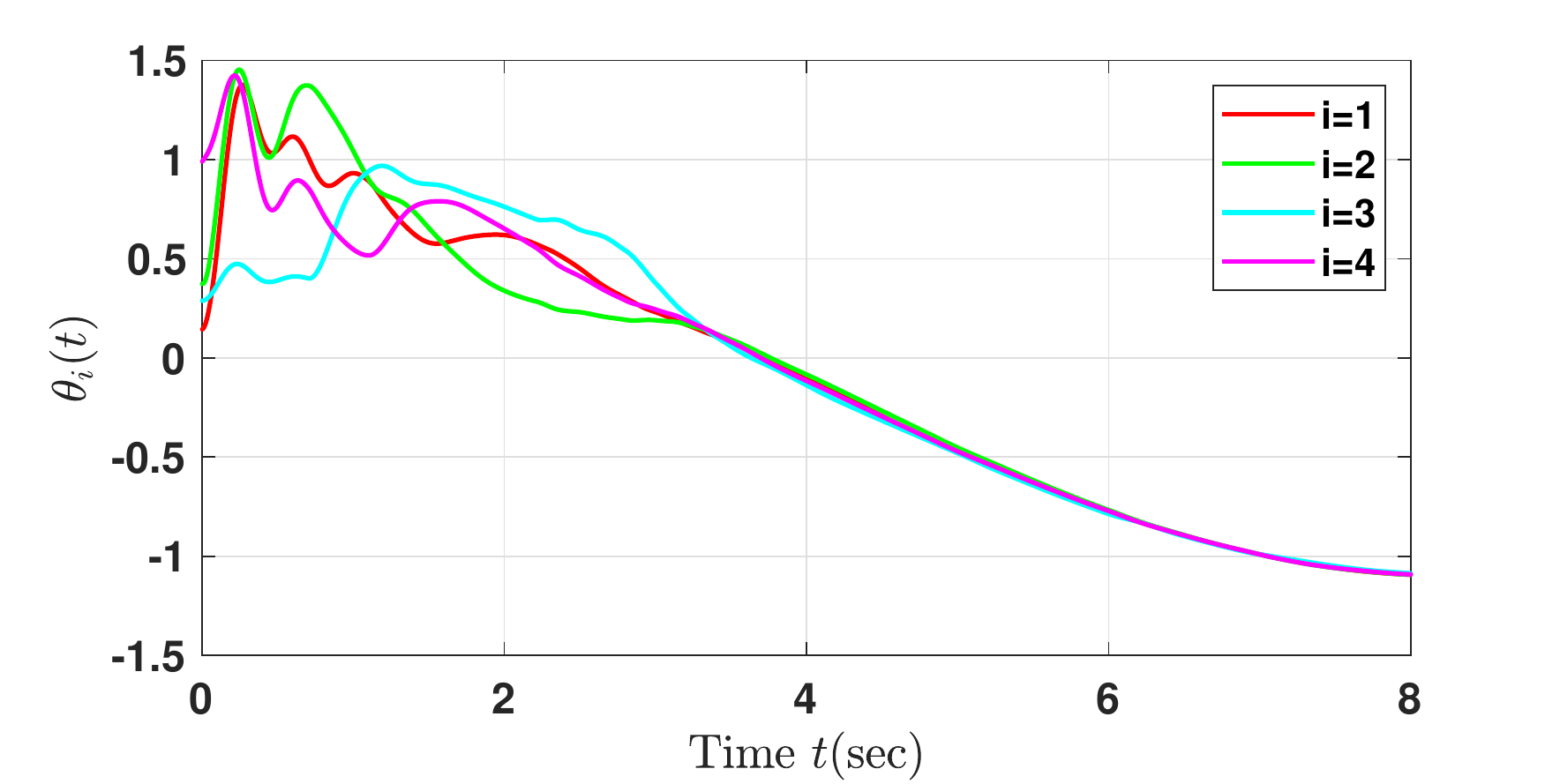}
\caption{The pitch angles of four agents under the reduced-order observer-based consensus protocol (\ref{RO_each_agent_controller}).}
\label{fig_theta_ro_trans}
\end{figure}

We show that both the full-state observer-based consensus protocol and the reduced-order observer-based consensus protocol enable four agents to achieve consensus when agents are perturbed by external disturbances in the form of square wave with period $2\pi$. Figs.~\ref{fig_phi_fs_trans}
and \ref{fig_theta_fs_trans} show the roll angles and pitch angles of these agents driven by the full-state observer-based consensus protocol with control gain matrices described by (\ref{FS_infinity_transient_control_gain}) and (\ref{FS_infinity_transient_observer_gain}). And Figs.~\ref{fig_phi_ro_trans} and \ref{fig_theta_ro_trans} present the corresponding trajectories of four agents driven by the reduced-order observer-based consensus protocol with control gain matrices given by (\ref{RO_infinity_transient_control_gain_barF}), (\ref{RO_infinity_transient_control_gain_G}), (\ref{RO_infinity_transient_control_gain}), (\ref{RO_infinity_transient_control_gain_T}), (\ref{RO_infinity_transient_control_gain_Q1}), and  (\ref{RO_infinity_transient_control_gain_Q2}).
It can be seen from these figures that all agents can reach consensus though they are perturbed by the external disturbances.

\begin{figure}
\includegraphics[width=4in]{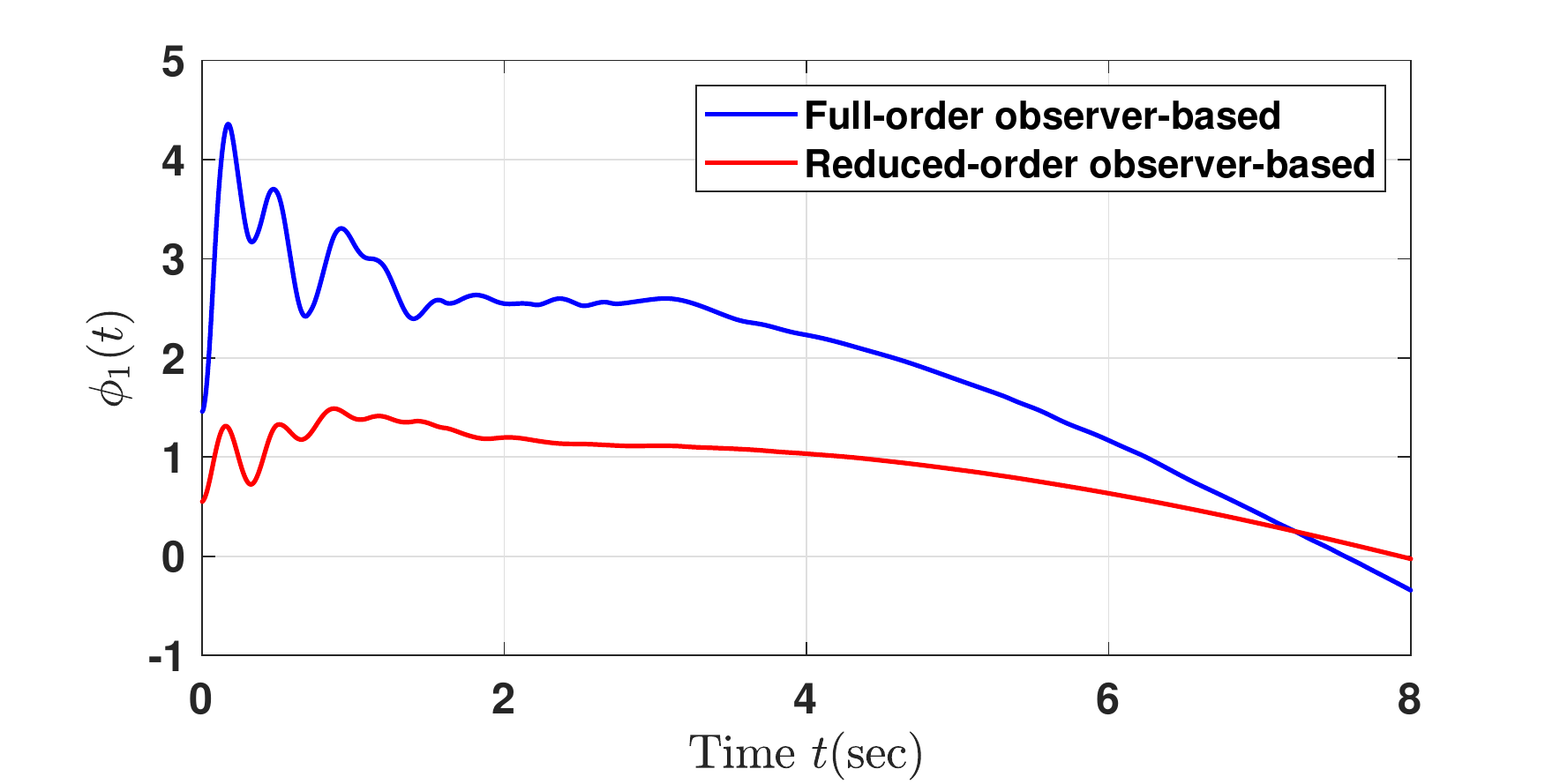}
\caption{The roll angles of the fist agent under the full-state observer-based consensus protocol (blue line) and the reduced-order observer-based control protocol (red line).}
\label{fig_COM_PHI1}
\end{figure}

\begin{figure}
\includegraphics[width=4in]{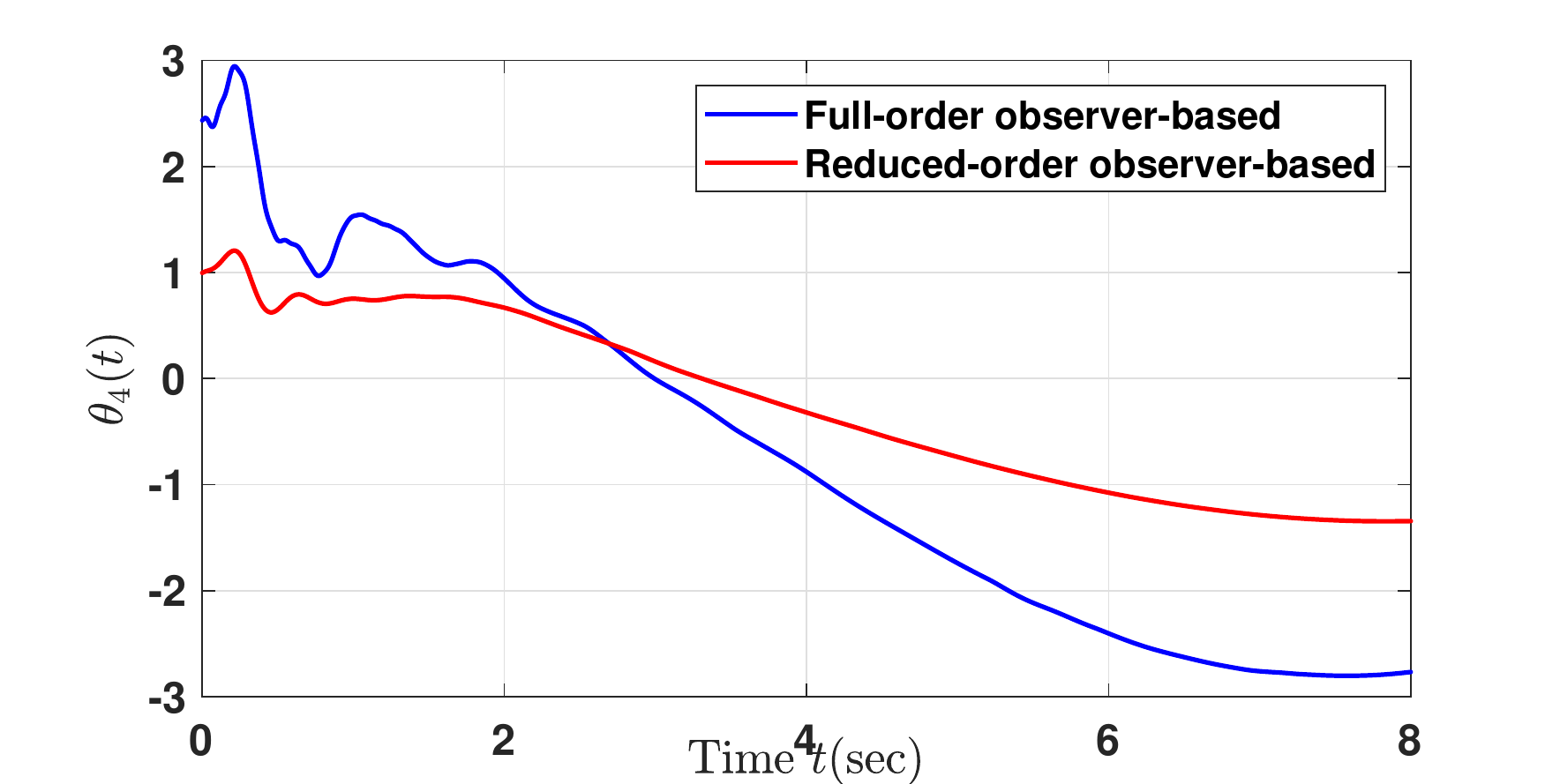}
\caption{The pitch angles of the forth agent under the full-state observer-based consensus protocol (blue line) and the reduced-order observer-based control protocol (red line).}
\label{fig_COM_THETA4}
\end{figure}

Now we make a comparison of the performance of grouped agents under the above two protocols. Fist, we are interested in the consensus time required by the grouped agents, which also reflects the convergence rate. It follows from Figs.~\ref{fig_phi_fs_trans}
and \ref{fig_theta_fs_trans} that agents need almost $4$ seconds to reach consensus driven by the full-state observer-based consensus protocol, while they cost almost $6$ seconds under the reduced-order observer-based consensus protocol, which can be seen from Figs.~\ref{fig_phi_ro_trans} and \ref{fig_theta_ro_trans}. Next, we move our attention to the transient performance of agents under both protocols. We draw the trajectories of the same state of some agent in one picture to achieve our goal. We choose the first agent's roll angle and the fourth agent's pitch angle to illustrate that reduced-order observer-based consensus protocol allow agents to have better transient behaviour, which are respectively shown in Figs.~\ref{fig_COM_PHI1} and \ref{fig_COM_THETA4}. It can be seen from Fig.~\ref{fig_COM_THETA4} that the overshoot of the forth agent's pitch angle under the full-state observer-based consensus protocol is almost three times larger than that under the reduced-order observer-based consensus protocol. What is more, the trajectory of the pitch angle under the full-state observer-based consensus protocol is much more oscillatory than that under the reduced-order observer-based consensus protocol. Precisely, the full-state observer-based consensus protocol makes the agent to eliminate the oscillation in almost $4$ seconds, while the reduced-order observer-based consensus protocol just costs the agent almost $1$ second. Figs.~\ref{fig_COM_PHI1}
shows that the first agent behaves similarly. Thus, we can conclude that compared with the full-state observer-based consensus protocol, the reduced-order observer-based consensus protocol brings better transient performance to agents, but converges in a relatively slower rate. By analysing the poles of the closed-loop system, we find that the reduced-order observer-based consensus protocol places the poles closer from the real axis in the left half s-plane, which results in less oscillation cycles. Therefore, we can conclude that both protocols guarantee robustness against external disturbance, while the reduced-order observer-based consensus protocol constructed here performs better in terms of the transient performance.

\section{Conclusions}\label{sec_conclusion}
This paper has considered the $H_{\infty}$ consensus control integrated with transient performance improvement problem for agents with general linear dynamics. Two classes of observer-based protocols, which have been proposed to deal with this problem, solely require the output information of nodes. Thus, they are immune to the absence of entire state information. Algorithms have been designed to construct these protocols, which are applicable in large-scale networks. At the end of this manuscript, numerical simulations have been presented to show the effectiveness of these algorithms. In the future work, we will investigate the robust stochastic consensus control problem of networked agents with uncertainties in their dynamics.








\begin{thebibliography}{1}



\bibitem{Liu2015finite}
Y. Liu, and Z. Geng, ``Finite-time formation control for linear multi-agent systems: A motion planning approach", \emph{Systems and Control Letters}, vol. 85, no. 11, pp. 54-60, 2015.

\bibitem{Zhou1996robust}
K. Zhou, J. C. Doyle, and K. Glover, Robust and Optimal Control, \emph{Prentice Hall}, 1996.








\bibitem{robust_H_infty_containment_Jia_2015}
P. Wang, and Y. Jia, ``Robust ${H}_{\infty}$ containment control for second-order multi-agent systems with nonlinear dynamics in directed networks", \emph{Nerocomputing}, vol. 153,  pp. 235--241, 2015.






\bibitem{wang2013h}
J. Wang, Z. Duan, Y. Zhao, G. Qin, and Y. Yan, ``${H}_{\infty}$ and ${H}_{2}$ control of multi-agent systems with transient performance improvement", \emph{International Journal of Control}, vol. 86, no. 12, pp. 2131--2145, 2013.


\bibitem{robust_consensus_Zeng_2015}
W. Huang, J. Zeng, and H. Sun, ``Robust consensus for linear multi-agent systems with mixed uncertainties", \emph{Systems and Control Letters}, vol. 76, pp. 56--65, 2015.

\bibitem{Hinf_Htwo_perfor_region_2010}
Z. Li, Z. Duan, and G. Chen, ``On ${H}_{\infty}$ and ${H}_{2}$ performance regions of multi-agent systems", \emph{Automatica}, vol. 57, no. 1, pp.
213--224, 2011.

\bibitem{Wang2015distributed}
J. Wang, Z. Duan, G. Wen, and G. Chen, ``Distributed robust control of uncertain linear multi-agent systems", \emph{International Journal of Robust and Nonlinear Control}, vol. 25, pp. 2162--2179, 2015.

\bibitem{Wang2015consensus}
Q. Wang, H. Di, Z. Duan, and J. Wang, ``Consensus tracking control with transient performance improvement for a group of unmanned aerial vehicles subject to faults and parameter uncertainty", \emph{International Journal of Control}, DOI: 10.1080/00207179.2017.1370555, 2017.


\bibitem{lin2017distributedconstrained}
P. Lin, and W. Ren, ``Distributed ${H}_{\infty}$ constrained consensus problem", \emph{Systems and Control Letters}, vol. 104, pp. 45-48.



\bibitem{zhao2012hinfinity}
Y. Zhao, Z. Duan, G. Wen, and G. Chen, ``Distributed ${H}_{\infty}$ consensus of multi-agent systems: A performance region-based approach", \emph{International Journal of Control}, vol. 85, no. 3, pp. 332--341, 2012.


\bibitem{Boyd2006randomized}
S. Boyd, A. Ghosh, B. Prabhakar, and D. Shah, ``Randomized gossip algorithms", \emph{{IEEE} Transactions on Information Theory}, vol. 52, no. 6, pp. 2508--2530, 2006.



\bibitem{Matei2009consensus_problems}
I. Matei, N. Martins, and J. S. Baras, ``Consensus problems with directed Markovian communication patterns", \emph{In Proceedings of the American Control Conference}, pp. 1298--1303, St. Louis, MO, USA.




\bibitem{Li2015containment}
W. Li, L. Xie, and J. Zhang, ``Containment control of leader-following multi-agent systems with Markovian switching network topologies and measurement noises", \emph{Automatica}, vol. 51, pp. 263--267, 2015.





\bibitem{Kan2016leader}
Z. Kan, J. M. Shea, and W. E. Dixon, ``Leader-follower containment control over directed random graphs", \emph{Automatica}, vol. 66, pp. 56-62.




\bibitem{Wang2015seeking}
Y. Wang, L. Cheng, W. Ren, Z. Hou, and M. Tan, ``Seeking consensus in networks of linear agents: Communication noises and Markovian switching topologies", \emph{{IEEE} Transactions on Automatic Control}, vol. 60, no. 5, pp. 1374--1379.


\bibitem{you2013consensus}
K. You, Z. Li, and L. Xie, ``Consensus condition for linear multi-agent systems over randomly switching topologies", \emph{Automatica}, vol. 49, no. 10, pp. 3125--3132, 2013.




\bibitem{Pan2017distributedcooperative}
Y. Pan, H. Werner, Z. Huang, and M. Bartels, ``Distributed cooperative control of leader-follower multi-agent systems under packet dropouts for quadcopters", \emph{Systems and Control Letters}, vol. 106, pp. 47--57, 2017.








\bibitem{Mu2014L}
X. Mu, B. Zheng, and K. Liu, ``${L}_{2}$-${L}_{\infty}$
containment control of multi-agent systems with Markovian switching topologies and non-uniform time-varying delays", \emph{{IET} Control Theory and Applications}, vol. 8, iss. 10, pp. 863--872, 2014.







\bibitem{Wang2017distributed}
J. Wang, Z. Duan, G. Wen, and Y. Hao, ``Distributed ${H}_{\infty}$ control of multi-agent systems over randomly switching topologies", \emph{Proceedings of the 2017 International Workshop on Computer Science and Networks (IWCSN)}, pp. 156-160, Doha, Qatar, December 2017.




\bibitem{IJC_2010}
D. V. Balandin, and M. M. Kogan, ``{LMI}-based ${H}_{\infty}$-optimal control with transients", \emph{International Journal of Control}, vol. 83, no. 8, pp. 1664--1673, 2010.



\bibitem{wang2018h}
J. Wang, Z. Duan, and G. Wen, ``${H}_{\infty}$ containment control for
multi-agent systems over Markovian switching topologies", \emph{Proceedings of the 33rd Youth Academic Annual Conference of Chinese Association of Automation (YAC 2018)}, Nanjing, China, May 2018.















\bibitem{report_TYAA_flying_vehicle_2005}
B. Chen, ``Fourth half yearly progress report for the {TYIA} 2003 project on nonlinear control methods for flight control systems of flying vehicles", \emph{DSTA, NUS, and ST Aerospace}, 2005.











































\end{thebibliography}
\end{document}